\documentclass[12pt]{article}
\usepackage{geometry}                
\usepackage{graphicx}
\usepackage{amsmath, amssymb, amsthm}
\usepackage{lscape}
\usepackage{multirow}
\usepackage{url}
\usepackage[numbers,square,sort]{natbib}
\usepackage[ruled, vlined]{algorithm2e}

\usepackage{setspace}
\newtheorem{lem}{Lemma}

\title{A unified modeling approach for the static-dynamic uncertainty strategy in stochastic lot-sizing}
\author{Roberto Rossi,${}^{*,1}$ Onur A. Kilic,${}^2$ S. Armagan Tarim,${}^2$\\
${}^1$Business School, University of Edinburgh, Edinburgh, UK\\
roberto.rossi@ed.ac.uk\\
${}^2$Institute of Population Studies, Hacettepe University, Ankara, Turkey\\
\{onuralp,armagan.tarim\}@hacettepe.edu.tr
}
\date{}                                           

\begin{document}
\maketitle

\begin{abstract}
In this paper, we develop mixed integer linear programming models to compute near-optimal policy parameters for the non-stationary stochastic lot sizing problem under Bookbinder and Tan's static-dynamic uncertainty strategy. Our models build on piecewise linear upper and lower bounds of the first order loss function. We discuss different formulations of the stochastic lot sizing problem, in which the quality of service is captured by means of backorder penalty costs, non-stockout probability, or fill rate constraints. These models can be easily adapted to operate in settings in which unmet demand is backordered or lost. The proposed approach has a number of advantages with respect to existing methods in the literature: it enables seamless modelling of different variants of the above problem, which have been previously tackled via ad-hoc solution methods; and it produces an accurate estimation of the expected total cost, expressed in terms of upper and lower bounds. Our computational study demonstrates the effectiveness and flexibility of our models.\\
{\bf keywords: } 
stochastic lot sizing; static-dynamic uncertainty; first order loss function; non-stockout probability; fill rate; penalty cost; piecewise linearisation\\
{\bf Corresponding author}: Roberto Rossi, University of Edinburgh Business School, EH8 9JS, Edinburgh, United Kingdom, 
phone: +44(0)131 6515239, email: roberto.rossi@ed.ac.uk
\end{abstract}

\newpage


\section{Introduction}

Consider the non-stationary stochastic lot sizing problem --- the stochastic extension of the well-known dynamic lot sizing problem \citep{Wagner1958}. This is a finite-horizon periodic review single-item single-stocking location inventory control problem in which demand is stochastic and non-stationary. \citet{bt88} discuss three main control strategies that can be adopted in stochastic lot sizing problem: static, static-dynamic, and dynamic uncertainty. The static uncertainty strategy is rather conservative. The decision maker determines both timing and size of orders at the very beginning of the planning horizon. A less conservative strategy is the static-dynamic uncertainty, in which inventory reviews are fixed at the beginning of the planning horizon, while associated order quantities are decided upon only when orders are issued. The dynamic uncertainty strategy allows the decision maker to decide dynamically at each time period whether or not to place an order and how much to order. Each of these strategies has different advantages and disadvantages. For instance, the dynamic uncertainty strategy is known to be cost-optimal \citep{Scarf1960}. The static uncertainty is appealing in material requirement planning systems, for which order synchronisation is a key concern \citep{citeulike:12394813}. The static-dynamic uncertainty strategy has advantages in organising joint replenishments and shipment consolidation \citep{spp98,Mutlu2010}. 

In this study, we focus our attention on the static-dynamic uncertainty strategy, which offers a stable replenishment plan while effectively hedging against uncertainty \citep{Kilic2011,Tunc2012}. An important question regarding the static-dynamic uncertainty strategy is how to determine order quantities at inventory review periods when a replenishment schedule is given. In this context, \citet{citeulike:10430648} showed that it is optimal to determine order quantities by means of an order-up-to policy. This result leads to the following characterisation of the static-dynamic uncertainty strategy: at each review period, the decision maker observes the actual inventory position (i.e. on-hand inventory, plus outstanding orders, minus backorder) and places an order so as to increase the inventory position up to a given order-up-to level. Key decisions for the static-dynamic uncertainty strategy include an inventory review schedule and an order-up-to level for each review period --- these decisions must be fixed at the beginning of the planning horizon. 

We build on recently introduced piecewise linear upper and lower bounds for the first order loss function and its complementary function \cite{citeulike:12518455}, which are based on distribution independent bounding techniques from stochastic programming: Jensen's and Edmundson-Madanski's inequalities \cite[][p. 167-168]{citeulike:695971}. In contrast to earlier works in the literature, we show that these bounds can be used to estimate inventory holding costs, backorder costs and/or service levels, and that they translate into readily available lower and upper bounds on the optimal expected total costs. Furthermore, for the special case in which demand is normally distributed, the model relies on standard linearisation parameters provided in \cite{citeulike:12518455}. 

Our analysis leads to a unified modelling approach that captures several variants of the problem and that is based on standard mixed-integer linear programming models. Some of these variants have been previously addressed in the literature, whereas some other have not. More specifically, we consider different assumptions on the way unsatisfied demand is modelled: backorder and lost sales. We also consider different service quality measures commonly employed in the inventory control literature \citep[see e.g.][pp. 244---246]{spp98}: penalty cost per unit short per period, non-stockout probability\footnote{a lower bound on the non-stockout probability in any period over the planning horizon} ($\alpha$ service level), cycle fill rate\footnote{a lower bound on the expected fraction of demand that is routinely satisfied from stock for each replenishment cycle} ($\beta^{\text{cyc}}$ service level), and fill rate\footnote{a lower bound on the expected fraction of demand that is routinely satisfied from stock over the planning horizon} ($\beta$ service level). 

Our contributions to the inventory control literature are the following:
\begin{itemize}
\item we develop enhanced MILP formulations that enable seamless modelling of the non-stationary stochastic lot sizing problem under each of the four measures of service quality discussed;
\item in contrast to other approaches in the literature our MILP models bound from above and below the cost of an optimal plan by using a piecewise linear approximation of the loss function; by increasing the number of segments, precision can be improved ad libitum;
\item we discuss how to build these MILP models for the case in which demand in each period follows a generic probability distribution;
\item for the special case in which demand in each period is normally distributed, we demonstrate how the MILP formulations can be conveniently constructed via standard linearisation coefficients;
\item we discuss for the first time in the literature how to handle the case in which demand that occurs when the system is out of stock is lost, i.e. lost sales;
\item we discuss the first MILP formulation in the literature that captures the case in which service quality is modelled using a standard $\beta$ service level in line with the definition found in many textbooks on inventory control, such as \citet[][]{citeulike:11028967,spp98,Axsater2006}.
\item we present an extensive computational study to show that (i) the optimality gap shrinks exponentially fast as the number of segments in the piecewise linearisation increases, (ii) the number of segments adopted only marginally affects computational efficiency.
\end{itemize}

\section{Literature survey}\label{sec:literature}
Due to its practical relevance, a large body of literature has emerged on the static-dynamic uncertainty strategy over the last decades.
Here, we review some key contributions which are of particular importance in the context of the current work, and reflect upon our contribution. 

To keep our discussion focused, since all the models we will discuss in the following sections operate under a static-dynamic uncertainty strategy, in our literature review we do not survey works related to the static uncertainty \cite[see e.g.][]{citeulike:12827714,citeulike:8923542,citeulike:8763804,citeulike:10379346} and to the dynamic uncertainty strategy, \cite[see e.g.][]{bollapragada1999}.

Early works on the stochastic lot sizing problem concentrated on easy-to-compute heuristics. \citet{sil78} and \citet{ask81} studied the problem under penalty costs, and proposed simple heuristics based on the least period cost method. These heuristics can be regarded as stochastic extensions of the well-known Silver-Meal heuristic \citep{Silver1973}. 

\citet{bt88} studied the problem under $\alpha$ service level constraints and introduced the terminology ``static uncertainty,'' ``dynamic uncertainty,'' and ``static-dynamic uncertainty.'' They developed a method that sequentially determines the timing of replenishments and corresponding order-up-to levels for the static-dynamic uncertainty strategy. Following this seminal work, a variety of further studies --- which significantly differ in terms of underlying service quality measures and modeling approaches --- aimed to determine the optimal replenishment schedule and order-up-to levels simultaneously under Bookbinder and Tan's static-dynamic uncertainty strategy. 

\citet{citeulike:12252354} discussed the first MILP formulation under $\alpha$ service level constraints. In contrast to \cite{bt88}, this formulation simultaneously determines the replenishment schedule and corresponding order-up-to levels. Efficient reformulations operating under the same assumptions were discussed in \cite{citeulike:9768278,citeulike:7156145}. \cite{citeulike:7156145} developed a state space augmentation approach; while \cite{citeulike:9768278} implemented a branch and bound algorithm. In addition, Constraint Programming reformulations based on a novel modelling tool, i.e. global chance constraints, were discussed in \cite{citeulike:6067172,citeulike:8772285}. Finally, an exact, although computationally intensive, Constraint Programming approach was discussed in \cite{citeulike:7875691}. Extensions to the case of a stochastic delivery lead time were discussed in \cite{citeulike:7288002,citeulike:12379660}.

\citet{citeulike:12317242} developed the first MILP formulation for the case in which service quality is modelled using a penalty cost scheme. \citet{citeulike:9668374} discussed an efficient Constraint Programming reformulation exploiting optimization oriented global stochastic constraints.

\citet{citeulike:10430648} discussed a dynamic programming solution algorithm and two ad-hoc heuristics named ``approximation'' and ``relaxation'' heuristics, respectively; the authors analyse both the penalty cost and the $\alpha$ service level cases. The ``approximation'' heuristic operates under the assumption that scenarios in which the actual stock exceeds the order-up-to-level for a given review are negligible and can be safely ignored; while the ``relaxation'' heuristics operates by relaxing those constraints in Tarim and Kingsman's model that force order sizes in each period to be nonnegative. 

\citet{citeulike:8766754} introduced an MILP formulation for the case in which service quality is modelled via $\beta^{\text{cyc}}$ service level constraints. 

A key issue in all the aforementioned studies is the computation of the true values of expected on-hand inventories and stock-outs, and thereby associated costs and/or service levels. These values can only be derived from the (complementary) first-order loss function of the demand \citep[see][p. 338]{citeulike:12819903} --- a non-linear function that cannot be readily embedded into the proposed MILP models. 

\citet{citeulike:12252354} --- but also \cite{citeulike:8772285, citeulike:9768278, citeulike:6067172, citeulike:7875691,citeulike:7288002,citeulike:12379660,citeulike:7156145} --- bypassed the issue by approximating the expected on-hand inventory by the expected inventory position. This approach could work well for inventory systems that operate under a very high non-stockout probability. However, it may result in highly sub-optimal solutions when the probability of observing a stock-out is not negligible. 

\citet{citeulike:12317242} used a piecewise linear approximation of the standard loss function to approximate expected holding and penalty costs. The piecewise linear function is fitted to the nonlinear cost function by using an approach that minimises the maximum absolute approximation error. The power of this approach is that the piecewise linearisation is based on standard linearisation coefficients that can be computed once and then reused for any normally distributed demand. Unfortunately, the approximation proposed may either over or underestimate the original cost therefore it becomes hard to assess how far a given solution may be from the true optimal one, i.e. its optimality gap. Furthermore, this approximation is not easily extended to demands following a generic distribution or to models operating under service level measures rather than a penalty cost scheme.

\citet{citeulike:10430648} approximations require ad-hoc algorithms and cannot be easily extended to handle $\beta$ service level constraints or lost sales.

\citet{citeulike:8766754} tabulates the complementary first order loss function and then uses binary variables to retrieve the holding cost associated with a given replenishment plan. A similar tabulation is employed to enforce the prescribed $\beta^{\text{cyc}}$ service level. However, this tabulation is carried out by considering each possible replenishment cycle\footnote{A replenishment cycle is the time interval between two successive inventory reviews.} independently. This strategy disregards cost and service level dependencies that may exist among successive replenishment cycles. For this very same reason, it cannot be employed to model classic $\beta$ service level constraints. To the best of our knowledge, no formulation exists in the literature for the case in which service quality is modelled via standard $\beta$ service level constraints \cite[][]{citeulike:11028967,spp98,Axsater2006}.

The issue of computing expected on-hand inventories and stock-outs is topical in inventory control, as witnessed by a number of recent works \citep[see e.g.][]{citeulike:7766622,citeulike:3169436}. All approaches surveyed above address particular instances of the problem. The contribution of this paper is unique and novel in the sense that it introduces a unified modelling approach for static-dynamic uncertainty strategy based on linear approximations of the first order loss function. Furthermore, this 
unified modelling approach can be used to address the issue of computing static-dynamic uncertainty policy parameters under lost sales, which has not been addressed yet in the literature.

\section{Piecewise linearisation of loss functions}\label{sec:piecewise}

For convenience, a list of all symbols used in the rest of the paper is provided in Appendix I. Consider a random variable $\omega$ with expected value $\tilde{\omega}$ and a scalar variable $x$. The first order loss function is defined as
$\mathcal{L}(x,\omega)=\mbox{E}[\max(\omega-x,0)]$,
where $\mbox{E}$ denotes the expected value. The complementary first order loss function is defined as
$\widehat{\mathcal{L}}(x,\omega)=\mbox{E}[\max(x-\omega,0)]$.
It is know that there is a close relationship between these two functions, as stated in the following lemma.
\begin{lem}[\cite{citeulike:12819903} p. 338, C.5]\label{thm:relationship_fol_com_1}
\begin{equation}
\mathcal{L}(x,\omega)=\widehat{\mathcal{L}}(x,\omega)-(x-\tilde{\omega}).
\end{equation}
\end{lem}
The first order loss function and its complementary function play a key role in inventory models, since they are essential to compute expected holding and penalty costs, as well as a number of service measure such as the $\beta^{\text{cyc}}$ and $\beta$ service levels.

A common approach in computing near-optimal control parameters of the static-dynamic uncertainty strategy is to formulate the problem as a certainty equivalent MILP. Unfortunately, the loss function is non-linear and cannot be easily embedded in MILP models. To overcome this issue, we adopt a piecewise linearisation approach similar to the one recently discussed in \citet{citeulike:12518455}. This approach is based on classical inequalities from stochastic programming: Jensen's and Edmundson-Madanski inequalities \cite[][p. 167-168]{citeulike:695971} and can be applied to random variables following a generic probability distribution.

It is known that both the first order loss function $\mathcal{L}(x,\omega)$ and its complementary function $\widehat{\mathcal{L}}(x,\omega)$ are convex in $x$ regardless of the distribution of $\omega$. For this reason, both Jensen's lower bound and Edmundson-Madanski upper bound are applicable to the first order loss function  and its complementary function. More formally, let $g_{\omega}(\cdot)$ denote the probability density function of $\omega$ and consider a partition of the support $\Omega$ of $\omega$ into $W$ disjoint compact subregions $\Omega_1,\ldots,\Omega_W$. We define, for all $i=1,\ldots,W$ 
\begin{equation}\label{eq:conditional_expectation}
p_i=\Pr\{\omega\in \Omega_i\}=\int_{\Omega_i} g_{\omega}(t)\,dt~~~\text{and}~~~
\mbox{E}[\omega|\Omega_i]=\frac{1}{p_i}\int_{\Omega_i} t g_{\omega}(t)\,dt
\end{equation}

\begin{lem}\label{lem:piecewise_linear_lb}
For the complementary first order loss function the lower bound $\widehat{\mathcal{L}}_{lb}(x,\omega)$, where
\[\widehat{\mathcal{L}}(x,\omega)\geq\widehat{\mathcal{L}}_{lb}(x,\omega)=\sum_{i=1}^W p_i \max(x-\mbox{E}[\omega|\Omega_i],0)\]
is a piecewise linear function with $W+1$ segments. The $i$-th linear segment of $\widehat{\mathcal{L}}_{lb}(x,\omega)$ is
\begin{equation}\label{eq:piecewise_linear_lb}
\widehat{\mathcal{L}}^i_{lb}(x,\omega)=x\sum_{k=1}^i p_k-\sum_{k=1}^i p_k \mbox{E}[\omega|\Omega_k]~~~~\mbox{E}[\omega|\Omega_{i}]\leq x\leq \mbox{E}[\omega|\Omega_{i+1}],
\end{equation}
where $i=1,\ldots,N$; furthermore, the $0$-th segment is $x=0$, $-\infty\leq x\leq \mbox{E}[\omega|\Omega_1]$.
\end{lem}
This lower bound is a direct application of Jensen's inequality. Let then $e_W$ denote the maximum approximation error for the lower bound in Lemma \ref{lem:piecewise_linear_lb} associated with a given partition comprising $W$ regions. A piecewise linear upper bound, i.e. Edmundson-Madanski's bound, can be obtained by shifting up the lower bound in Lemma \ref{lem:piecewise_linear_lb} by a value $e_W$.
\begin{lem}\label{lem:piecewise_linear_ub}
For the complementary first order loss function the upper bound $\widehat{\mathcal{L}}_{ub}(x,\omega)$, where
\[\widehat{\mathcal{L}}(x,\omega)\leq\widehat{\mathcal{L}}_{ub}(x,\omega)=\sum_{i=1}^W p_i \max(x-\mbox{E}[\omega|\Omega_i],0)+e_W\]
is a piecewise linear function with $W+1$ segments. The $i$-th linear segment of $\widehat{\mathcal{L}}_{ub}(x,\omega)$ is
\[\widehat{\mathcal{L}}^i_{ub}(x,\omega)=x\sum_{k=1}^i p_k-\sum_{k=1}^i p_k \mbox{E}[\omega|\Omega_k]+e_W~~~~\mbox{E}[\omega|\Omega_{i}]\leq x\leq \mbox{E}[\omega|\Omega_{i+1}],\]
where $i=1,\ldots,N$; furthermore, the $0$-th segment is $x=e_W$, $-\infty\leq x\leq \mbox{E}[\omega|\Omega_1]$.
\end{lem}

Having established these two results, we must then decide how to partition the support $\Omega$ in order to obtain good bounds. A number of works discussed how to obtain an optimal partitioning of the support under a framework that minimises the maximum approximation error \cite{citeulike:12820831,Imamoto2008}. In short, these works demonstrate that, in order to minimise the maximum approximation error, one must find parameters ensuring approximation errors at piecewise function breakpoints are all equal. This result unfortunately does not hold when optimal linearisation parameters must be found for a set of random variables. 

Consider a set of random variables $\omega_1,\ldots,\omega_n,\ldots,\omega_N$ and associated complementary first order loss functions $\widehat{\mathcal{L}}(x,\omega_1),\ldots, \widehat{\mathcal{L}}(x,\omega_N)$. 
From \eqref{eq:conditional_expectation} it is clear that, once all $p_i$ have been fixed, all $\mbox{E}[\omega_n|\Omega_k]$ are uniquely determined. The particular structure of \eqref{eq:piecewise_linear_lb} makes it, in principle, possible to compute standard $p_i$ coefficients for the whole set of random variables and then select the $\mbox{E}[\omega_n|\Omega_k]$ for a specific $\widehat{\mathcal{L}}_{lb}(x,\omega_n)$ via a binary selector variable $y_n$, that is
\begin{equation}\label{eq:piecewise_linear_lb_selector}
\widehat{\mathcal{L}}^i_{lb}(x,\omega)=x\sum_{k=1}^i p_k-\sum_{k=1}^i p_k \mbox{E}[\omega_n|\Omega_k] y_n~~~~
\begin{array}{l}
\mbox{E}[\omega_n|\Omega_{i}]\leq x\leq \mbox{E}[\omega_n|\Omega_{i+1}]\\ 
1\leq n\leq N
\end{array}
\end{equation}
where $\sum_{n=1}^N y_n=1$. These expressions generalise those discussed in \cite{citeulike:12518455}, which only hold for normally distributed random variables, and they are particularly useful, as we will see in Section \ref{sec:reformulations}, if one wants to develop an MILP model for computing Bookbinder and Tan's static-dynamic uncertainty policy parameters. 

Unfortunately, computing probability masses $p_1,\ldots,p_W$ that minimise the maximum approximation error over a set of random variables is a challenging task. As shown in \cite{rossi14} this is not a problem of convex optimisation as the one faced while computing optimal linearisation parameters for a single loss function. In this work, we will adopt two approximate strategies to compute good probability masses $p_1,\ldots,p_W$: a simple strategy that splits the support into $W$ regions with uniform probability mass; and a more refined local search strategy that uses a combination of simple random sampling and coordinate descent. Implementation details of these heuristic approaches are discussed in Appendix II.

\section{Enhanced MILP reformulations}\label{sec:reformulations}

In this section we demonstrate how the results presented so far can be used to derive enhanced a mixed integer programming formulations for the stochastic lot sizing problem under static-dynamic uncertainty strategy. First, we introduce the original stochastic programming formulation of the problem (Section \ref{sec:stoch_ls}). We then present MILP models for the problem under $\alpha$ service level constraints (Section \ref{sec:alpha_service}), a penalty cost scheme (Section \ref{sec:penalty}),  $\beta^{\text{cyc}}$ (Section \ref{sec:beta}) and $\beta$ (Section \ref{sec:beta_opt}) service level constraints. Finally, we discuss how to extend these models to a lost sales setting. For a complete overview of the models presented the reader may refer to Appendix III.

\subsection{Stochastic lot-sizing}\label{sec:stoch_ls}

The stochastic programming formulation of the non-stationary stochastic lot-sizing problem
was originally presented in \citep[][pp. 1097--1098]{bt88}. The formal problem definition is as follows. Customer demand $d_t$ in each period $t=1,\ldots,N$ is a random variable with probability density function $g_t(\cdot)$ and cumulative distribution function $G_t(\cdot)$. There are fixed and variable replenishment costs: the fixed cost is $a$ per order; the variable cost is $v$ per unit ordered. Negative orders are not allowed. A holding cost of $h$ is paid of each unit of inventory carried from one period to the next.  $I_0$ denotes the initial inventory level. Delivery lead-time is not incorporated in the model. When a stockout occurs, all demand is backordered and filled as soon as an adequate supply arrives. There is a service level constraint enforcing a non-stockout probability of at least $\alpha$ in each period --- this is known in the inventory control literature as ``$\alpha$ service level'' constraint \cite{spp98}. They finally also assume that the service level is set to a high value, i.e. $\alpha>0.9$ in order to incorporate management's perception of the cost of backorders, so that shortage costs can be safely ignored. The objective is to minimise the expected total cost $\mbox{E}[\mbox{TC}]$, which comprises fixed/variable ordering and holding costs. The resulting model is presented in Fig. \ref{fig:sp_model}.
\begin{figure}[h]
\begin{align}
\mbox{E}[\mbox{TC}]&=\min \int_{d_1}\int_{d_2}\ldots\int_{d_N} \sum_{t=1}^N (a\delta_t+h \max(I_t,0)+vQ_t)\times\\
&g_1(d_1)g_2(d_2)\ldots g_N(d_N)\,d(d_1)d(d_2)\ldots d(d_N)\nonumber\\
\nonumber\\
\mbox{subject to, }&\mbox{for }t=1,\ldots N\nonumber\\
\nonumber\\
&I_t=I_0+\sum_{i=1}^t(Q_i-d_i)\label{eq:inventory_conservation}\\
&\delta_t=\left\{
\begin{array}{ll}
1&\mbox{if }Q_t>0,\\
0&\mbox{otherwise}
\end{array}\right.\label{eq:reordering}\\
&\Pr\{I_t\geq0\}\geq\alpha\label{eq:service}\\
&Q_i\geq0,~\delta_t\in\{0,1\}
\end{align}
\caption{Stochastic programming formulation of the non-stationary stochastic lot-sizing problem.}
\label{fig:sp_model}
\end{figure}
In this model $I_t$ represents the inventory level at the end of a period; $\delta_t$ takes value 1 if an order is placed in period $t$; and $Q_t$ represents the order quantity in period $t$. Constraints \eqref{eq:inventory_conservation} are the inventory conservation constraints: inventory level at the end of period $t$ must be equal to the initial inventory $I_0$, plus all order received, minus all demand realised up to period $t$, since we assume --- in line with the original model in the literature --- that inventory cannot be disposed or returned to the supplier; constraints \eqref{eq:reordering} set $\delta_t$ to one if an order is placed in period $t$; finally, \eqref{eq:service} enforce the prescribed service level in each period. 

The above model can be easily modified to accommodate a penalty cost scheme, in place of the original $\alpha$ service level constraints, \citep[see e.g.][]{citeulike:12317242}. All one has to do is to drop \eqref{eq:service} and to replace the original objective function with the following one
\begin{align}
\mbox{E}[\mbox{TC}]&=\min \int_{d_1}\int_{d_2}\ldots\int_{d_N} \sum_{t=1}^N (a\delta_t+h \max(I_t,0)+b \max(-I_t,0)+vQ_t)\times\\
&g_1(d_1)g_2(d_2)\ldots g_N(d_N)\,d(d_1)d(d_2)\ldots d(d_N)\nonumber
\end{align}
where $b$ denotes the penalty cost per unit short per period. $\beta^{\text{cyc}}$ and $\beta$ service level formulations are obtained by replacing \eqref{eq:service} with an alternative service measure.

\subsection{$\alpha$ service level constraints}\label{sec:alpha_service}

We now consider the mixed integer programming formulation of \citet{citeulike:12252354} for computing near-optimal inventory control policy parameters under Bookbinder and Tan's static-dynamic uncertainty strategy. According to this strategy, inventory review times as well as their respective order-up-to-levels must be all fixed at the beginning of the planning horizon. However, actual order quantities are determined only after demand has been observed. 

In what follows, $M$ denotes a very large number and $\tilde{x}$ denotes the expected value of $x$. Tarim and Kingsman's model is presented in Fig. \ref{fig:milp_model}. 
\begin{figure}[h]
\begin{align}
\mbox{E}[\mbox{TC}]=-vI_0+v\sum_{t=1}^N \tilde{d}_t+&\min \sum_{t=1}^N (a\delta_t+h \tilde{I}_t)+v\tilde{I}_N\\
\nonumber\\
\mbox{subject to, }&\mbox{for }t=1,\ldots N\nonumber\\
\nonumber\\
&\tilde{I}_t+\tilde{d}_t-\tilde{I}_{t-1}\geq0\label{eq:inventory_conservation_milp}\\
&\tilde{I}_t+\tilde{d}_t-\tilde{I}_{t-1}\leq \delta_t M\label{eq:reordering_milp}\\
&\tilde{I}_t\geq \sum_{j=1}^t\left(G^{-1}_{d_{jt}}(\alpha)-\sum_{k=j}^t \tilde{d}_k\right)P_{jt}\label{eq:service_1}\\
&\sum_{j=1}^t P_{jt}=1\label{eq:service_2}\\
&P_{jt}\geq\delta_j-\sum_{k=j+1}^t\delta_k&j=1,\ldots,t\label{eq:service_3}\\
&P_{jt}\in\{0,1\}&j=1,\ldots,t\label{eq:binary_p_jt}\\
&\delta_t\in\{0,1\}\label{eq:binary_delta_t}
\end{align}
\caption{MILP formulation of the non-stationary stochastic lot-sizing problem under the static-dynamic uncertainty strategy \cite{citeulike:12252354}.}
\label{fig:milp_model}
\end{figure}
This certainty equivalent model comprises three sets of decision variables: $\tilde{I}_t$, representing the expected closing inventory level at the end of period $t$; $\delta_t$, a binary variable representing the inventory review decision at period $t$; and $P_{jt}$, a binary variable which is set to one if and only if the most recent inventory review before period $t$ was carried out in period $j$. By observing that, for a period $t$ in which an order is placed (i.e. $\delta_t=1$) the order-up-to-level $S_t$ is simply $S_t=\tilde{I}_t+\tilde{d}_t$, it follows that by solving the above model policy parameters are immediately obtained.

Constraints in the certainty equivalent model neatly reflect those in the original stochastic programming model. More specifically, \eqref{eq:inventory_conservation_milp} enforces the inventory conservation constraints; \eqref{eq:reordering_milp} is the reordering condition; and \eqref{eq:service_1}, \eqref{eq:service_2}, and \eqref{eq:service_3} enforce the prescribed service level $\alpha$. In \eqref{eq:service_1}, $G^{-1}_{d_{jt}}(\alpha)$ denotes the $\alpha$-quantile of the inverse cumulative distribution function of the random variable $d_j+\ldots+d_t$. Finally, the objective function is obtained by observing
\begin{equation}\label{eq:unit_cost}
\mbox{E}[v\sum_{t=1}^N Q_t]=-vI_0+v\sum_{t=1}^N \tilde{d}_t+v\tilde{I}_N,
\end{equation}
where $-vI_0+v\sum_{t=1}^N \tilde{d}_t$ is a constant. 

Following an assumption originally introduced by \citet{bt88}, \citet{citeulike:12252354} approximate the holding cost component in the original objective function, which we recall was $\mbox{E}[\max(I_t,0)]$, via the expression $h \tilde{I}_t$. To overcome this limitation of the model, we introduce two new sets of decision variables: $\tilde{I}^{lb}_t\geq0$ and $\tilde{I}^{ub}_t\geq0$ for $t=1,\ldots,N$, which represent, respectively, a lower and an upper bound to the true value of $\mbox{E}[\max(I_t,0)]$. The objective function then can be rewritten as
\begin{equation}\label{eq:obj_lb}
\mbox{E}[\mbox{TC}]=-vI_0+v\sum_{t=1}^N \tilde{d}_t+\min \sum_{t=1}^N (a\delta_t+h \tilde{I}^{lb}_t)+v\tilde{I}_N
\end{equation}
if our aim is to compute a lower bound for the cost of an optimal plan, or as
\begin{equation}\label{eq:obj_ub}
\mbox{E}[\mbox{TC}]=-vI_0+v\sum_{t=1}^N \tilde{d}_t+\min \sum_{t=1}^N (a\delta_t+h \tilde{I}^{ub}_t)+v\tilde{I}_N
\end{equation}
if our aim is to compute an upper bound for the cost of an optimal plan. We next discuss how to constrain $\tilde{I}^{lb}_t$ and $\tilde{I}^{ub}_t$. Our discussion applies to demand $d_k$, for $k=1,\ldots,N$ following a generic distribution. Consider the random variable $d_{jt}$ representing the convolution $d_j+\ldots+d_t$. We select a priori a number $W$ of adjacent regions $\Omega_i$ into which the support of $d_{jt}$ must be partitioned. As discussed, this partitioning will produce a piecewise linear approximation comprising $W+1$ segments. We also fix a priory the probability mass $p_i=\Pr\{d_{jt}\in\Omega_i\}$ that must be associated with each region $\Omega_i$. As discussed in Section \ref{sec:piecewise}, there are several possible strategies to assign probability masses $p_i$ to regions. For instance, we may ensure uniform coverage, i.e. all region must have the same probability mass, or we may select regions --- by using a heuristic or exact approach --- in order to minimise the maximum approximation error over all possible convolutions $d_{jt}$, for $t=1,\ldots,N$ and $j=1,\ldots,t$. Regardless of the strategy we adopt, once all $p_i$ are known, regions $\Omega_i$ are uniquely determined and the associated conditional expectation $\mbox{E}[d_{jt}|\Omega_i]$ can be immediately computed using off-the-shelf software. Finally, the maximum approximation error $e_W^{jt}$ associated with the linearisation of $\widehat{\mathcal{L}}(x,d_{jt})$ can be found by checking the linearisation error at the $W$ possible breakpoints of the piecewise linear function obtained. 
Having precomputed all these values for $t=1,\ldots,N$ and $j=1,\ldots,t$, we introduce the following constraints in the model
\begin{equation}\label{eq:piecewise_compl_loss_lb}
\tilde{I}^{lb}_t\geq (\tilde{I}_t+\sum_{j=1}^t\tilde{d}_{jt}P_{jt}) \sum_{k=1}^i p_k - \sum_{j=1}^t \left( \sum_{k=1}^i p_k\mbox{E}[d_{jt}|\Omega_i]\right) P_{jt}~~~t=1,\ldots,N;~i=1,\ldots,W 
\end{equation}
This expression follows from Lemma \ref{lem:piecewise_linear_lb} and closely resembles \eqref{eq:piecewise_linear_lb_selector}. Consider a replenishment in period $j$ covering periods $j,\ldots,t$ with associated order-up-to-level $S$. Our aim is to enforce $\tilde{I}^{lb}_t\geq \widehat{\mathcal{L}}^i_{lb}\left(S,d_{jt}\right)$ for $i=1,\ldots,W$, since $\tilde{I}^{lb}_t$ represents a lower bound for the expected positive inventory at the end of period $t$. By observing that $S=\tilde{I}_t+\tilde{d}_{jt}$, we obtain the above expression.
We then obtain $\tilde{I}^{ub}_t$  from $\tilde{I}_t$, by noting that a piecewise linear upper bound can be derived by adding the maximum estimation error to the Jensen's piecewise linear lower bound \cite{citeulike:12518455}.
\begin{equation}\label{eq:piecewise_compl_loss_ub}
\tilde{I}^{ub}_t\geq (\tilde{I}_t+\sum_{j=1}^t\tilde{d}_{jt}P_{jt}) \sum_{k=1}^i p_k + \sum_{j=1}^t \left(e_W^{jt} - \sum_{k=1}^i p_k\mbox{E}[d_{jt}|\Omega_i]\right)P_{jt}~~~
\begin{array}{l}
t=1,\ldots,N,\\
i=1,\ldots,W; 
\end{array}
\end{equation}
where $\tilde{I}^{ub}_t\geq\sum_{j=1}^t e_W^{jt}P_{jt}$ for $t=1,\ldots,N$. The special case in which demand in each period follows a normal distribution is discussed in Appendix IV.

\subsection{Penalty cost scheme}\label{sec:penalty}

The model discussed in Section \ref{sec:alpha_service} can be easily modified to accommodate a penalty cost $b$ per unit short per period in place of the $\alpha$ service level constraints discussed in \citet{citeulike:12252354}. This revised model resembles the one discussed in \citet{citeulike:12317242}. However, as discussed in the previous section, our formulation is more accurate, because the expected total cost of a plan can be now bounded from above and below. It is also more general, since the discussion in \citet{citeulike:12317242} is limited to normally distributed demand.

In the new model, we introduce two new sets of variables $\tilde{B}^{lb}_t\geq0$ and $\tilde{B}^{ub}_t\geq0$ for $t=1,\ldots,N$, which represent a lower and upper bound, respectively, for the true value of $\mbox{E}[-\min(I_t,0)]$ and thus allow us to compute lower and upper bounds for the expected backorders in each period. The objective function then becomes
\begin{equation}\label{eq:obj_pen_lb}
\mbox{E}[\mbox{TC}]=-vI_0+v\sum_{t=1}^N \tilde{d}_t+\min \sum_{t=1}^N (a\delta_t+h \tilde{I}^{lb}_t+b \tilde{B}^{lb}_t)+v\tilde{I}_N
\end{equation}
if our aim is to compute a lower bound for the cost of an optimal plan, or
\begin{equation}\label{eq:obj_pen_ub}
\mbox{E}[\mbox{TC}]=-vI_0+v\sum_{t=1}^N \tilde{d}_t+\min \sum_{t=1}^N (a\delta_t+h \tilde{I}^{ub}_t+b \tilde{B}^{ub}_t)+v\tilde{I}_N
\end{equation}
if our aim is to compute an upper bound for the cost of an optimal plan. Finally, we must remove constraints \eqref{eq:service_1}, since we are operating under a penalty cost scheme and not under a service level constraints. 

Once more, we assume that demand in each period follows a generic distribution; we obtain $\tilde{B}^{lb}_t$ and $\tilde{B}^{ub}_t$  from $\tilde{I}_t$ by exploiting Lemma \ref{thm:relationship_fol_com_1}.
\begin{equation}\label{eq:piecewise_loss_lb}
\tilde{B}^{lb}_t\geq - \tilde{I}_t + (\tilde{I}_t+\sum_{j=1}^t\tilde{d}_{jt}P_{jt}) \sum_{k=1}^i p_k - \sum_{j=1}^t \left( \sum_{k=1}^i p_k\mbox{E}[d_{jt}|\Omega_i]\right) P_{jt}~~~
\begin{array}{l}
t=1,\ldots,N,\\
i=1,\ldots,W; 
\end{array}
\end{equation}
where $\tilde{B}^{ub}_t\geq  - \tilde{I}_t$ and 
\begin{equation}\label{eq:piecewise_loss_ub}
\tilde{B}^{ub}_t\geq - \tilde{I}_t + (\tilde{I}_t+\sum_{j=1}^t\tilde{d}_{jt}P_{jt}) \sum_{k=1}^i p_k + \sum_{j=1}^t \left(e_W^{jt} - \sum_{k=1}^i p_k\mbox{E}[d_{jt}|\Omega_i]\right)P_{jt}~~~
\begin{array}{l}
t=1,\ldots,N,\\
i=1,\ldots,W; 
\end{array}
\end{equation}
where  $\tilde{B}^{ub}_t\geq  - \tilde{I}_t + \sum_{j=1}^t e_W^{jt} P_{jt}$. The case in which demand in each period follows a normal distribution is discussed in Appendix IV.

\subsection{$\beta^{\text{cyc}}$ service level constraints}\label{sec:beta}

The model discussed in Section \ref{sec:penalty} can be modified to accommodate $\beta^{\text{cyc}}$ service level constraints, defined in \cite{citeulike:8766754} as
\begin{equation}\label{beta_tempelmeier}
1 - \max_{i=1,\ldots,m}\left[\mbox{E}\left\{\frac{\text{Total backorders in replenishment cycle $i$}}{\text{Total demand in replenishment cycle $i$}}\right\}\right].
\end{equation}
$\beta^{\text{cyc}}$ represents a lower bound on the expected fraction of demand that is routinely satisfied from stock for each replenishment cycle. The revised model resembles the one discussed in \citet{citeulike:8766754}. However, we aim to stress that in this latter work the author enforces the prescribed service level by precomputing order-up-to-levels and cycle holding costs in a table, rather than using a piecewise linearisation of the loss function as we do. One of the advantages of our approach is that it is able to account for dependencies among opening stock levels and service levels of consecutive cycles. As we will show in Section \ref{sec:beta_opt}, the ability to capture these dependencies is important if we aim to generalise the model to a more classical definition of ``fill rate''.

The discussion below is distribution independent; we modify the model discussed in Section \ref{sec:penalty} by introducing service level constraints 
\begin{equation}\label{eq:fill_rate_lb}
\tilde{B}^{lb}_t\leq(1-\beta^{\text{cyc}})\sum_{j=1}^t P_{jt}\tilde{d}_{jt} ~~~t=1,\ldots,N,
\end{equation}
if our aim is to compute a lower bound for the cost of an optimal plan; or
\begin{equation}\label{eq:fill_rate_ub}
\tilde{B}^{ub}_t\leq(1-\beta^{\text{cyc}})\sum_{j=1}^t P_{jt}\tilde{d}_{jt} ~~~t=1,\ldots,N,
\end{equation}
if our aim is to compute an upper bound for the cost of an optimal plan. Constraints \eqref{eq:fill_rate_lb} and \eqref{eq:fill_rate_ub} directly follow from \eqref{beta_tempelmeier}.
Finally, the objective function is \eqref{eq:obj_lb} if our aim is to compute a lower bound for the cost of an optimal plan, or \eqref{eq:obj_ub} if our aim is to compute an upper bound. 

\subsection{$\beta$ service level constraints}\label{sec:beta_opt}

The model discussed in \citet{citeulike:8766754} captures a definition of fill rate that is not conventional in the inventory literature. This issue has been discussed in \citet{citeulike:12518457}. To date, no modelling strategy exists for the conventional fill rate under a static-dynamic uncertainty control policy. In this section, we introduce an alternative MILP reformulation that captures a definition of $\beta$ service level that is in line with the definition found in many textbooks on inventory control \citet[][]{citeulike:11028967,spp98,Axsater2006}. 

In \citet{Axsater2006}, the author defines $\beta$ service level as ``the expected fraction of demand satisfied immediately from stock on hand''. In the context of finite horizon inventory models, e.g. \citet{Chen2003,Thomas2005} this definition is formalized as
\begin{equation}\label{beta}
1 - \mbox{E}\left\{\frac{\text{Total backorders within the planning horizon}}{\text{Total demand within the planning horizon}}\right\},
\end{equation}
The static-dynamic uncertainty strategy divides the finite planning horizon into a number, say $m$, of consecutive replenishment cycles. We can re-write (\ref{beta}) by taking these into account as
\begin{equation}\label{beta_cycle}
1 - \mbox{E}\left\{\frac{\sum_{i=1}^m\text{Total backorders within the $i$'th replenishment cycle}}{\sum_{i=1}^m\text{Total demand within the $i$'th replenishment cycle}}\right\}.
\end{equation}
However, in \cite{citeulike:8766754} the same $\beta$ service level is imposed on each and every cycle within the planning horizon. Clearly \eqref{beta_cycle} is different from \eqref{beta_tempelmeier}: the original definition imposes a $\beta$ service level throughout the planning horizon, whereas the definition in \cite{citeulike:8766754} imposes a $\beta$ service level on each replenishment cycle within the planning horizon independently. This latter definition is thus more restrictive. Adopting the original definition may therefore bring substantial cost benefits. To illustrate this, one may consider instances discussed in \cite[][p. 191]{citeulike:8766754}, for which the cost reduction with respect to policies obtained via the model discussed in Section \ref{sec:beta} ranges from 0.4\% to 6.4\%.

We modify as follows the model in Section \ref{sec:beta} to implement the classical measure outlined in \eqref{beta_cycle}. Also in this cases the changes discussed below are distribution independent. We introduce two new set of nonnegative variables $\tilde{C}^{lb}_t$ and $\tilde{C}^{ub}_t$ for $t=0,\ldots,N$. These variables express a lower and an upper bound, respectively, to the expected total backorders within the replenishment cycle that ends at period $t$, if any exists. Hence, $\tilde{C}^{lb}_t$ (resp. $\tilde{C}^{ub}_t$) should be equal to $\tilde{B}^{lb}_t$ (resp. $\tilde{B}^{ub}_t$), if $t$ is the last period of a replenishment cycle;  otherwise $\tilde{C}^{lb}_t$ (resp. $\tilde{C}^{ub}_t$) should be equal to 0. We enforce this fact as follows:
\begin{equation}\label{eq:end_of_cycle_backorders_lb}
\tilde{C}^{lb}_t \geq \tilde{B}^{lb}_t  - (1-\delta_{t+1})\sum_{k=1}^t \tilde{d}_t~~~t=0,\ldots,N-1,
\end{equation}
\begin{equation}\label{eq:end_of_cycle_backorders_ub}
\tilde{C}^{ub}_t \geq \tilde{B}^{ub}_t  - (1-\delta_{t+1})\sum_{k=1}^t \tilde{d}_t ~~~t=0,\ldots,N-1.
\end{equation}
where $\tilde{B}^{lb}_0=\tilde{B}^{ub}_0=\tilde{C}^{lb}_0=\tilde{C}^{ub}_0=I_0$.
Finally, we must ensure that $\tilde{C}^{lb}_N=\tilde{B}^{lb}_N$ and $\tilde{C}^{ub}_N=\tilde{B}^{ub}_N$. We then use these new variables to build constraint
\begin{equation}\label{eq:classical_fill_rate_lb}
\sum_{t=1}^N\tilde{C}^{lb}_t\leq(1-\beta)\sum_{t=1}^N \tilde{d}_t
\end{equation}
which will replace \eqref{eq:fill_rate_lb}, if our aim is to compute a lower bound for the cost of an optimal plan; and constraint
\begin{equation}\label{eq:classical_fill_rate_ub}
\sum_{t=1}^N\tilde{C}^{ub}_t\leq(1-\beta)\sum_{t=1}^N \tilde{d}_t
\end{equation}
which will replace \eqref{eq:fill_rate_ub}, if our aim is to compute an upper bound. Constraints \ref{eq:classical_fill_rate_lb} and \ref{eq:classical_fill_rate_ub} directly follow from \ref{beta_cycle}.

\section{Lost sales}\label{sec:ls}

In this section we briefly sketch the extension of the models discussed in the previous section to the case in which demand that occurs when the system is out of stock is lost, i.e. lost sales setting. The discussion is purportedly short, since the derived models are quite similar to those already presented and there are only few adjustments that are necessary to adapt our models to this new settings.

Under lost sales we need to take into account the fact that if inventory drops to zero, demand that occurs until the next order arrives will not be met. Under a static-dynamic uncertainty control policy, this means the actual order quantity will never exceed the order-up-to-level. 

In this settings it is crucial to set up the model in such a way as to account for the opportunity cost associated with units of demand that are not met by a given control policy. For this reason, we must introduce a new parameter $s$ that represents the selling price of a product; we then let $m=s-v$ be the margin --- i.e. unit selling price minus unit ordering cost --- for an item sold. The resulting stochastic programming model under lost sales is shown in Fig. \ref{fig:sp_model_lost_sales}.
\begin{figure}[h]
\begin{align}
\mbox{E}[\mbox{TP}]=sI_0+&\max\int_{d_1}\int_{d_2}\ldots\int_{d_N} \left[\sum_{t=1}^N (mQ_t-a\delta_t-h \max(I_t,0))-s\max(I_N,0)\right]\times\\
&g_1(d_1)g_2(d_2)\ldots g_N(d_N)\,d(d_1)d(d_2)\ldots d(d_N)\nonumber\\
\nonumber\\
\mbox{subject to, }&\mbox{for }t=1,\ldots N\nonumber\\
\nonumber\\
&I_t+d_t-I_{t-1}\geq0\label{eq:inventory_conservation_ls}\\
&Q_t= I_t+d_t-\max(I_{t-1},0)\label{eq:order_qty_ls}\\
&\delta_t=\left\{
\begin{array}{ll}
1&\mbox{if }Q_t>0,\\
0&\mbox{otherwise}
\end{array}\right.\label{eq:reordering_ls}\\
&\Pr\{I_t\geq0\}\geq\alpha\label{eq:service_ls}\\
&Q_t\geq0,~\delta_t\in\{0,1\}
\end{align}
\caption{Stochastic programming formulation of the non-stationary stochastic lot-sizing problem under lost sales.}
\label{fig:sp_model_lost_sales}
\end{figure}
In this model $\mbox{E}[\mbox{TP}]$ represents the expected total profit, which we aim to maximise. $\delta_t$ is a binary decision variable that is set to one if we order items at period $t$, i.e. constraints \eqref{eq:reordering_ls}. $Q_t$ represents the order quantity in period $t$, which must be greater or equal to zero. $I_t$ is a random variable that represents the inventory level at the end of period $t$; despite lost sales, for convenience, we assume that $I_t$ may take negative values. 
In the objective function, we multiply the margin $m$ by the number of items $Q_t$ ordered in period $t$; we then subtract ordering cost $a$, if an order is placed in period $t$ (i.e. $\delta_t=1$), and the holding cost $h$ on items that remain in stock at the end of period $t$. There is a further term $-s\max(I_N,0)$ to reflect the fact that items in stock at the end of the planning horizon will not be sold and thus the associated selling price $s$ should not be included in the total profit. Term $\max(I_{t-1},0)$ in constraints \eqref{eq:order_qty_ls} makes sure that the order quantity does not include any lost sale from the previous period. Constraints \eqref{eq:inventory_conservation_ls} ensure that the inventory level at the end of period $t$ is greater or equal to the inventory level at the end of period $t-1$ plus the realised demand in period $t$; this makes sure that items in stock in a period and not sold are brought to the next period. Constraints \eqref{eq:service_ls} enforce an $\alpha$ service level and can be easily replaced by other service measures such as a $\beta^{\text{cyc}}$ or a $\beta$ service level. A formulation under a penalty cost scheme is also easily obtained by removing the service level constraints and by modifying the objective function as illustrated in previous sections.

An MILP formulation of the problem in Fig. \ref{fig:sp_model_lost_sales} under the static-dynamic uncertainty strategy is shown in Fig. \ref{fig:milp_model_ls}.
\begin{figure}[h!]
\begin{align}
\mbox{E}[\mbox{TP}]=&sI_0+\max \sum_{t=1}^N (m\tilde{Q}_t - a\delta_t-h \tilde{I}^{lb}_t)-s\tilde{I}^{lb}_N\label{eq_obj_ls_lb}\\
\nonumber\\
\mbox{subject to, }&\mbox{for }t=1,\ldots N\nonumber\\
\nonumber\\
&\tilde{I}_t+\tilde{d}_t-\tilde{I}_{t-1}\geq0\label{eq:inventory_conservation_milp_ls}\\
&\tilde{Q}_t\geq \tilde{I}_t+\tilde{d}_t-\tilde{I}^{lb}_{t-1} + (1-\delta_t) M\label{eq:reordering_milp_ls_1a}\\
&\tilde{Q}_t\leq \tilde{I}_t+\tilde{d}_t-\tilde{I}^{lb}_{t-1} - (1-\delta_t) M\label{eq:reordering_milp_ls_1b}\\
&\tilde{I}_t+\tilde{d}_t-\tilde{I}_{t-1} \leq \tilde{Q}_t M\label{eq:reordering_milp_ls_3}\\
&\tilde{Q}_t\leq \delta_t M\label{eq:reordering_milp_ls_2}\\
&\tilde{I}_t\geq \sum_{j=1}^t\left(G^{-1}_{d_{jt}}(\alpha)-\sum_{k=j}^t \tilde{d}_k\right)P_{jt}\label{eq:service_1_ls}\\
&\sum_{j=1}^t P_{jt}=1\label{eq:service_2_ls}\\
&P_{jt}\geq\delta_j-\sum_{k=j+1}^t\delta_k&j=1,\ldots,t\label{eq:service_3_ls}\\
&\mbox{constraints \eqref{eq:piecewise_compl_loss_lb}}\nonumber\\
&P_{jt}\in\{0,1\}&j=1,\ldots,t\label{eq:binary_p_jt_ls}\\
&\tilde{Q}_t\geq0,~\delta_t\in\{0,1\}\label{eq:binary_delta_t_ls}
\end{align}
\caption{MILP formulation of the non-stationary stochastic lot-sizing problem under the static-dynamic uncertainty strategy and lost sales.}
\label{fig:milp_model_ls}
\end{figure}
In this model, by taking expectations, constraints \eqref{eq:inventory_conservation_ls} translate into \eqref{eq:inventory_conservation_milp_ls}. Constraints \eqref{eq:order_qty_ls} translate into \eqref{eq:reordering_milp_ls_1a} and \eqref{eq:reordering_milp_ls_1b}.  Terms $\max(I_{t},0)$ in the objective function and in constraints \eqref{eq:order_qty_ls} can be handled by using an auxiliary variable $\tilde{I}^{lb}_t$ that represents a lower bound for $\mbox{E}[\max(I_{t},0)]$ computed as before via a piecewise linearisation of the complementary first order loss function, i.e. constraints \eqref{eq:piecewise_compl_loss_lb}. Constraints \eqref{eq:reordering_milp_ls_3} ensures that the expected inventory level at the end of period $t$ is greater than the expected inventory level at the end of period $t-1$ plus the expected demand in period $t$ if and only if an order has been placed in period $t$, i.e. $\tilde{Q}_t>0$. Constraints \eqref{eq:reordering_ls} translate into \eqref{eq:reordering_milp_ls_2}. Finally, service level constraints \eqref{eq:service_ls} translate into \eqref{eq:service_1_ls}, \eqref{eq:service_2_ls} and \eqref{eq:service_3_ls} following a strategy similar to the one in \cite{citeulike:12252354}, which we illustrated in Section \ref{sec:alpha_service}. The model presented can be used to compute an upper bound for $\mbox{E}[\mbox{TP}]$ --- note that underestimating buffer stocks, i.e. $\tilde{I}^{lb}_t$ leads to lower holding costs and to an overestimation of the expected order quantity and associated margins $m\tilde{Q}_t$ in the objective function. If we aim to compute a lower bound instead, all occurrences of $\tilde{I}^{lb}_t$ should be replaced by $\tilde{I}^{ub}_t$ and constraints \eqref{eq:piecewise_compl_loss_lb} should be replaced by constraints \eqref{eq:piecewise_compl_loss_ub}. Other MILP formulations under $\beta^{\text{cyc}}$ and $\beta$ service levels are obtained in a similar fashion, since only the service level constraints of the model are affected by this change. A penalty cost formulation is also easily obtained, by dropping the service level constraints and by modifying the objective function as follows 
\begin{equation}\label{eq_obj_ls_pen_lb}
\mbox{E}[\mbox{TP}]=sI_0+\max \sum_{t=1}^N (m\tilde{Q}_t - a\delta_t-h \tilde{I}^{lb}_t-b \tilde{B}^{lb}_t)-s\tilde{I}^{lb}_N
\end{equation}
A detailed overview of these models is given in Appendix III. Finally, a remark that should be made is that the model presented in Section \ref{sec:penalty} charges penalty cost on a ``per unit short per unit time'' basis. This might be not appropriate under lost sales, since it is common practice to charge penalty cost on a ``per unit short'' basis under this setting. This problem can be easily overcome by charging the penalty cost $b$ not on $\tilde{B}_t^{lb}$, which bounds the expected units short at the end of each period, but on $\tilde{C}_t^{lb}$, which bounds the expected units short at the end of each replenishment cycle --- as defined in Section \ref{sec:beta_opt}.

\section{Computational experience}\label{sec:comp_exp}

In this section we present an extensive computational analysis of the models previously discussed. The experiments below were conducted by using CPLEX 12.3 on a 2.13 Ghz Intel Core 2 Duo with 4GB of RAM.

The aim of our computational analysis is twofold. First, we investigate the behaviour of the optimality gap for all models presented when the number of segments adopted for the piecewise linear approximation of the loss function increases. The term optimality gap is used here to denote the difference between the upper and lower bounds for the expected total cost obtained via Edmundson-Madanski and Jensen's bounds, respectively. Second, we investigate the computational efficiency of our models and how the number of segments adopted in the piecewise linear approximation of the loss function impacts solution times. 

\subsection{Test bed}

\begin{figure}[p!]
\centering
\includegraphics[type=eps,ext=.eps,read=.eps,width=0.85\columnwidth]{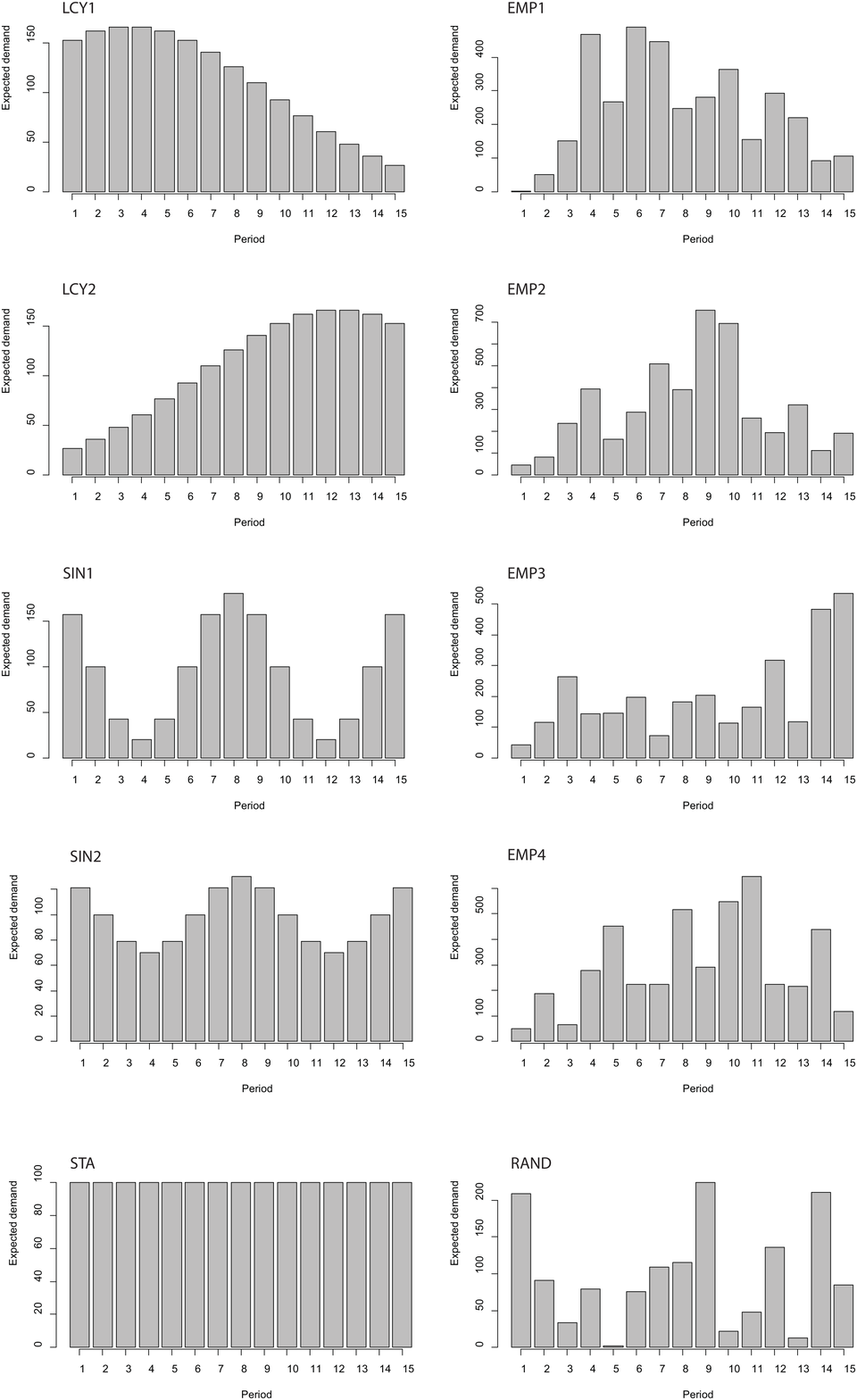}
\caption{Demand patterns in our computational analysis; the values presented denote the expected demand $\tilde{d}_t$ in each period $t$ of the planning horizon.}
\label{fig:cap_1}
\end{figure}

We consider a test bed comprising 810 instances. More specifically, we carried out a full factorial analysis under the following factors. We considered ten different demand patterns illustrated in Fig. \ref{fig:cap_1}. The patterns include two life cycle patterns (LCY1 and LCY2), two sinusoidal patterns (SIN1 and SIN2), stationary (STA) and random (RAND) patterns, and four empirical patterns derived from demand data in \cite{doi:10.1287/opre.44.1.131} (EMP1,\ldots,EMP4). 
Fixed ordering cost $a$ takes values in $\{500,1000,2000\}$; while proportional unit cost $v$ takes values in $\{2,5,10\}$. For models under service level constraints, we let the prescribed service level range in $\{0.8,0.9,0.95\}$; for model under a penalty cost scheme we let the penalty cost range in $\{2,5,10\}$.

\subsection{Backorders}

In this section, we concentrate on models presented in Section \ref{sec:stoch_ls}. Recall that in these models, when a stockout occurs, all demand is backordered and filled as soon as an adequate supply arrives. We first investigate the behaviour of the optimality gap and of the solution time for normally distributed demand. We then extend the analysis to the case in which demand in different periods follow different probability distributions.

\subsubsection{Normal distribution}
In this section, we assume that demand $d_t$ in each period $t$ to be normally distributed with mean $\tilde{d}_t$  and standard deviation $\sigma_{d_t}$, and let the coefficient of variation $c=\tilde{d}_t/\sigma_{d_t}$; we let $c\in\{0.10,0.20,0.30\}$. Expected values of the demand in each period are illustrated in Fig. \ref{fig:cap_1} for each of the ten patterns considered. As discussed, when demand is normally distributed, general purpose linearisation parameters can be precomputed and immediately used in our models \cite{citeulike:12518455}. 

\begin{figure}[p!]
\centering
\includegraphics[type=eps,ext=.eps,read=.eps,width=1\columnwidth]{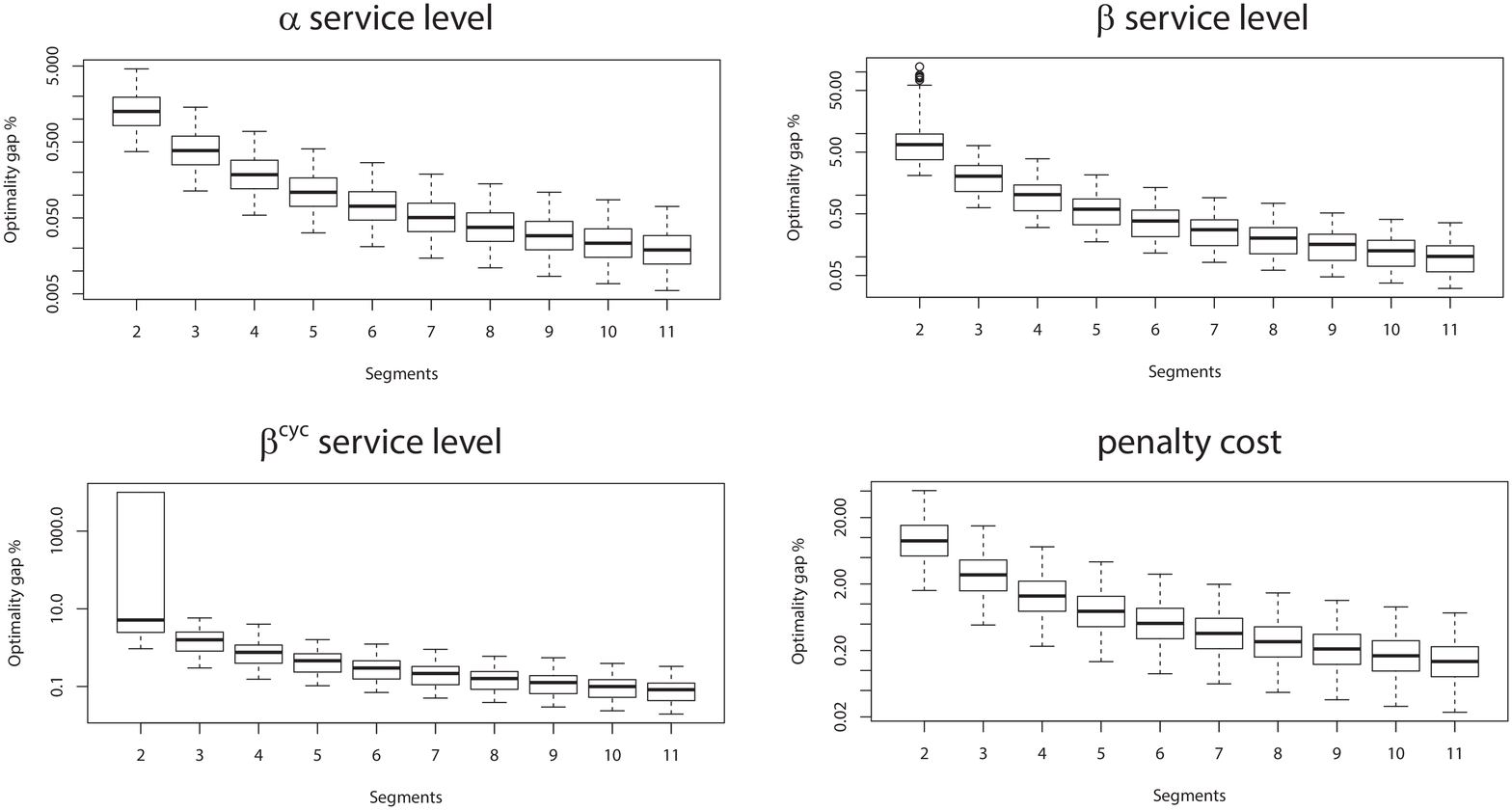}
\caption{Boxplots illustrating the optimality gap trend for different number of segments used in the piecewise linear approximation.}
\label{fig:cap_2}
\end{figure}

In Fig. \ref{fig:cap_2} we report, for each model discussed, boxplots illustrating the optimality gap trend for different number of segments used in the piecewise linear approximation. 
It should be noted that the y-axis is displayed in logarithmic scale. This shows that the optimality gap shrinks exponentially fast in the number of segments regardless of the model or parameter setting considered. 

It is interesting to observe that a number of instances were found infeasible by the MILP model under $\beta^{\text{cyc}}$ service level when two segments for the piecewise approximation were used --- note the very large optimality gap. This is due to the fact that with only two segments the approximation error for the $\beta^{\text{cyc}}$ service level was too large and no order-up-to-level could be found to enforce a service level as high as the prescribed one. 

\begin{figure}[p!]
\centering
\includegraphics[type=eps,ext=.eps,read=.eps,width=1\columnwidth]{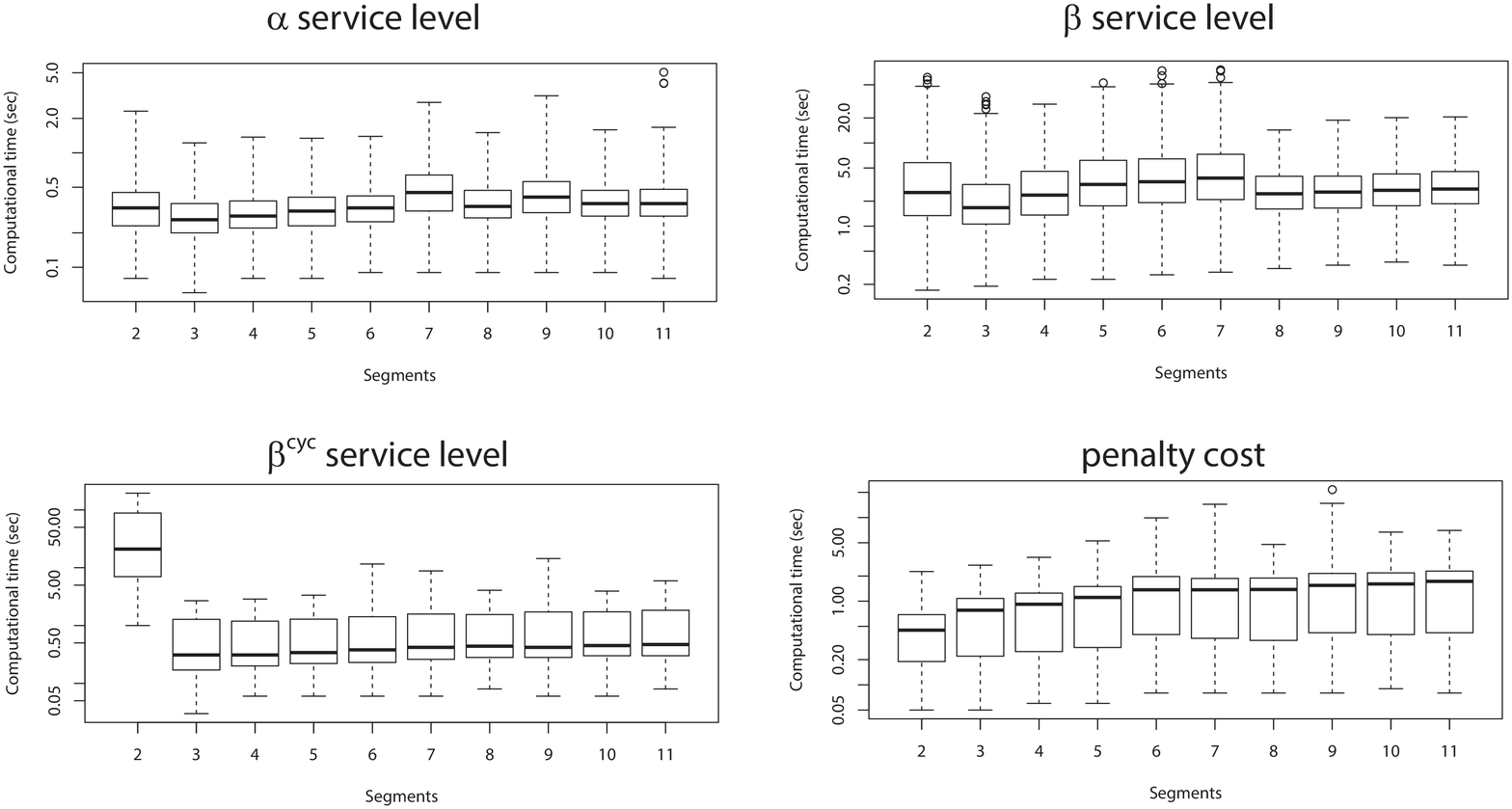}
\caption{Boxplots illustrating the computational time trend for different number of segments used in the piecewise linear approximation. }
\label{fig:cap_3}
\end{figure}

In Fig. \ref{fig:cap_3} we report, for each model discussed, boxplots illustrating the computational time trend for different number of segments used in the piecewise linear approximation. 
Computational times are only slightly affected by the number of segments in the approximation. Furthermore, all instances could be solved in just a few seconds. 

\subsubsection{Generic distribution}

We now extend the analysis to the case in which demand in different periods follow different probability distributions. The test bed remains the same previously analysed. However, the probability distribution of the demand in each period is now detailed in Table \ref{tab:prob_dist_demand}. To the best of our knowledge, this is the first study in which a numerical analysis on demand that is not normally distributed is presented.
\begin{table}
\begin{tabular}{ll}
Period	&Distribution\\
\hline
1,5,9,13	&normal with mean $\tilde{d}_t$ and $\sigma_{d_t}=0.3\cdot\tilde{d}_t$\\
2,6,10,14	&Poisson with mean $\tilde{d}_t$\\
3,7,11,15	&exponential with mean $\tilde{d}_t$\\
4,8,12	&uniformly distributed in $[0,2\tilde{d}_t]$
\end{tabular}
\caption{Probability distribution of the demand in each period of the planning horizon.}
\label{tab:prob_dist_demand}
\end{table}
We limit our analysis to the model that operates under penalty cost scheme (Section \ref{sec:penalty}). We do this for two reasons: firstly, because our local search procedure is  computationally quite intensive and a larger test bed would have taken considerable time; secondly, because the model that operates under penalty cost scheme embeds both the first order loss function --- employed to compute expected shortages --- and the complementary first order loss function --- employed to compute holding costs; therefore we expect that results obtained for this model to be sufficiently representative of the overall degree of approximation attained by our piecewise linear approximation of these two functions.

\begin{figure}[h!]
\centering
\includegraphics[type=eps,ext=.eps,read=.eps,width=1\columnwidth]{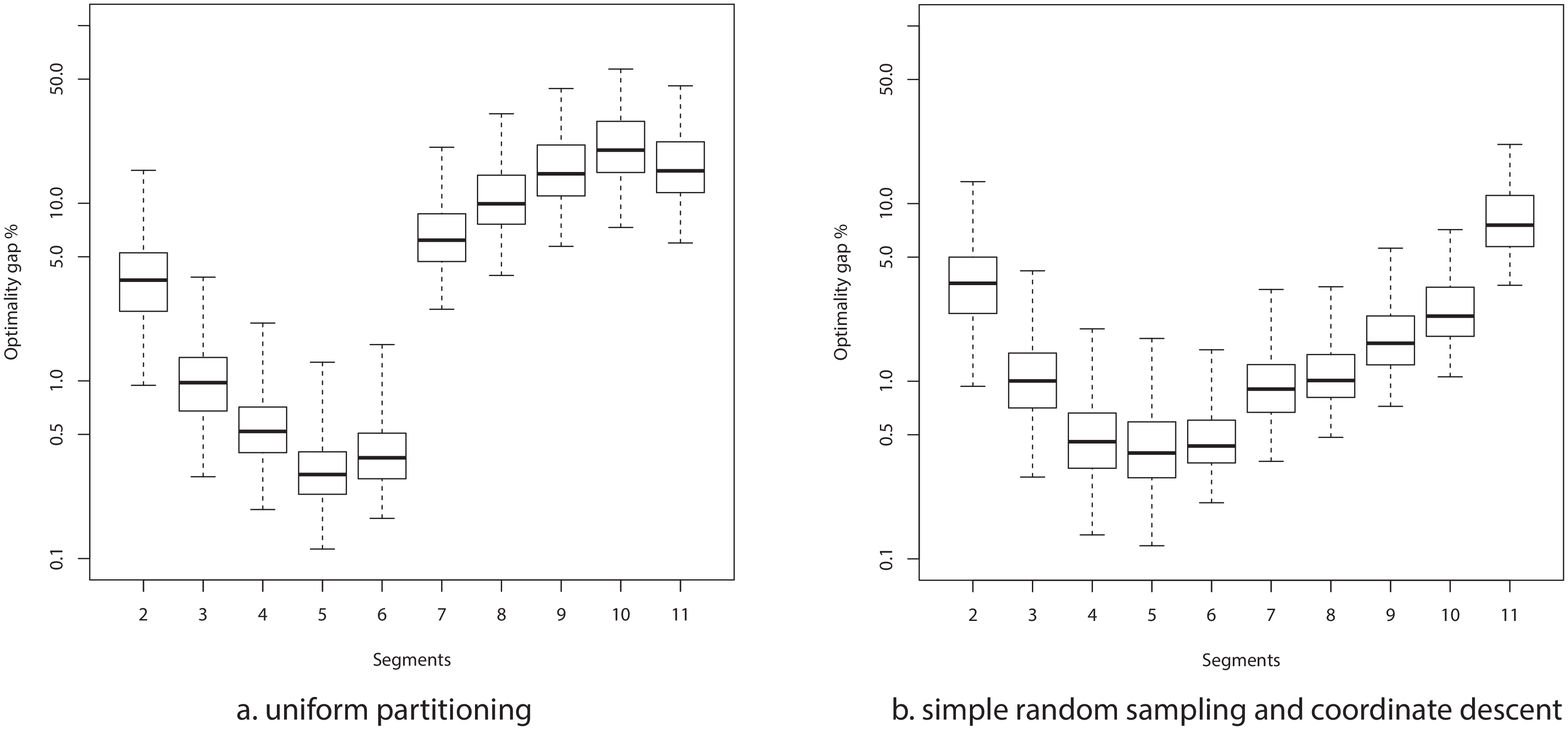}
\caption{Boxplots illustrating the optimality gap trend for different number of segments used in the piecewise linear approximation.}
\label{fig:cap_4}
\end{figure}

As discussed in Section \ref{sec:piecewise} when demand follows a generic probability distribution we must compute dedicated linearisation parameters for the piecewise first order loss function. In Fig. \ref{fig:cap_4} we consider two possible strategies for computing these parameters. (Fig. \ref{fig:cap_4} - a) illustrates results for the case in which we split the support of the demand $\omega$ uniformly into $W$ disjoint compact subregions $\Omega_1,\ldots,\Omega_W$ such that $p_i=\Pr\{\omega\in \Omega_i\}=1/W$. (Fig. \ref{fig:cap_4} - b) illustrates results for the case in which a local search procedure involving a coordinate descent from the most promising partition obtained via simple random sampling is employed to find a good partition into $W$ disjoint compact subregions that minimises the maximum approximation error. 

The results presented reveal that a uniform partitioning of the support of the demand provides acceptable results --- optimality gap around 0.5\% when the number of segments in the linearisation is low. However, as the number of segments increases, the performance this strategy deteriorates exponentially --- note that axis are in logarithmic scale. Local search may help attaining better performances when the number of segments increases. In our specific example, our simple local search strategy, which used a population size $\mathcal{S}=500W$ for the simple random sampling and a step size $\epsilon=0.002$ in the coordinate descent, attained up to an order of magnitude improvement in the optimality gap.  However, as discussed in \cite{rossi14}, this problem is computationally challenging and future research should investigate global optimisation algorithms to compute an optimal partitioning.

A final remark should be made on the computational efficiency of our approach under a generic demand distribution. If one adopts a uniform partitioning of the support then the computational efficiency of the model is essentially identical to the case in which demand is normally distributed. In other words, most of the instances can be solved in few seconds as shown in Fig. \ref{fig:cap_3}. If, however, a local search procedure is used to partition the support of the demand then the computational performance of our approach will of course depend on the efficiency of this local search strategy. In our analysis, the local search strategy effectively became the bottleneck, since finding good linearisation parameters for a specific instance took a time that varied from few minutes up to an hour for large $W$. Also in this case we believe that substantial improvements in computational performances may be achieved via more effective global optimisation algorithms.

\subsection{Lost sales}

We finally extended our analysis to the case in which demand that occurs when the system is out of stock is lost. Models that operate under this settings were discussed in Section \ref{sec:ls}. We analyse the case in which demand $d_t$ in each period $t$ is normally distributed with coefficient of variation $c\in\{0.10,0.20,0.30\}$. 

\begin{figure}[p!]
\centering
\includegraphics[type=eps,ext=.eps,read=.eps,width=1\columnwidth]{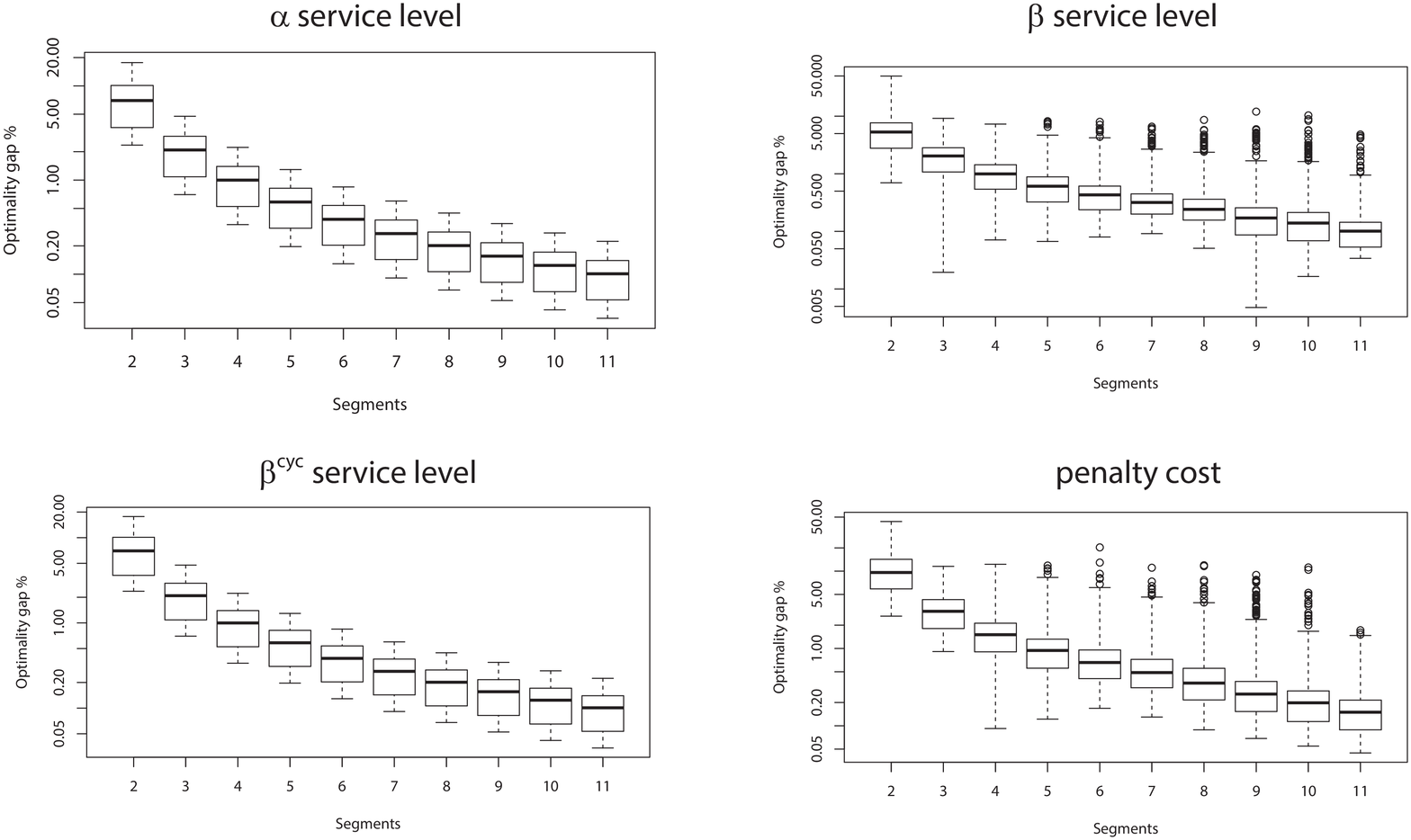}
\caption{Boxplots illustrating the optimality gap trend for different number of segments used in the piecewise linear approximation.}
\label{fig:cap_5}
\end{figure}

\begin{figure}[p!]
\centering
\includegraphics[type=eps,ext=.eps,read=.eps,width=1\columnwidth]{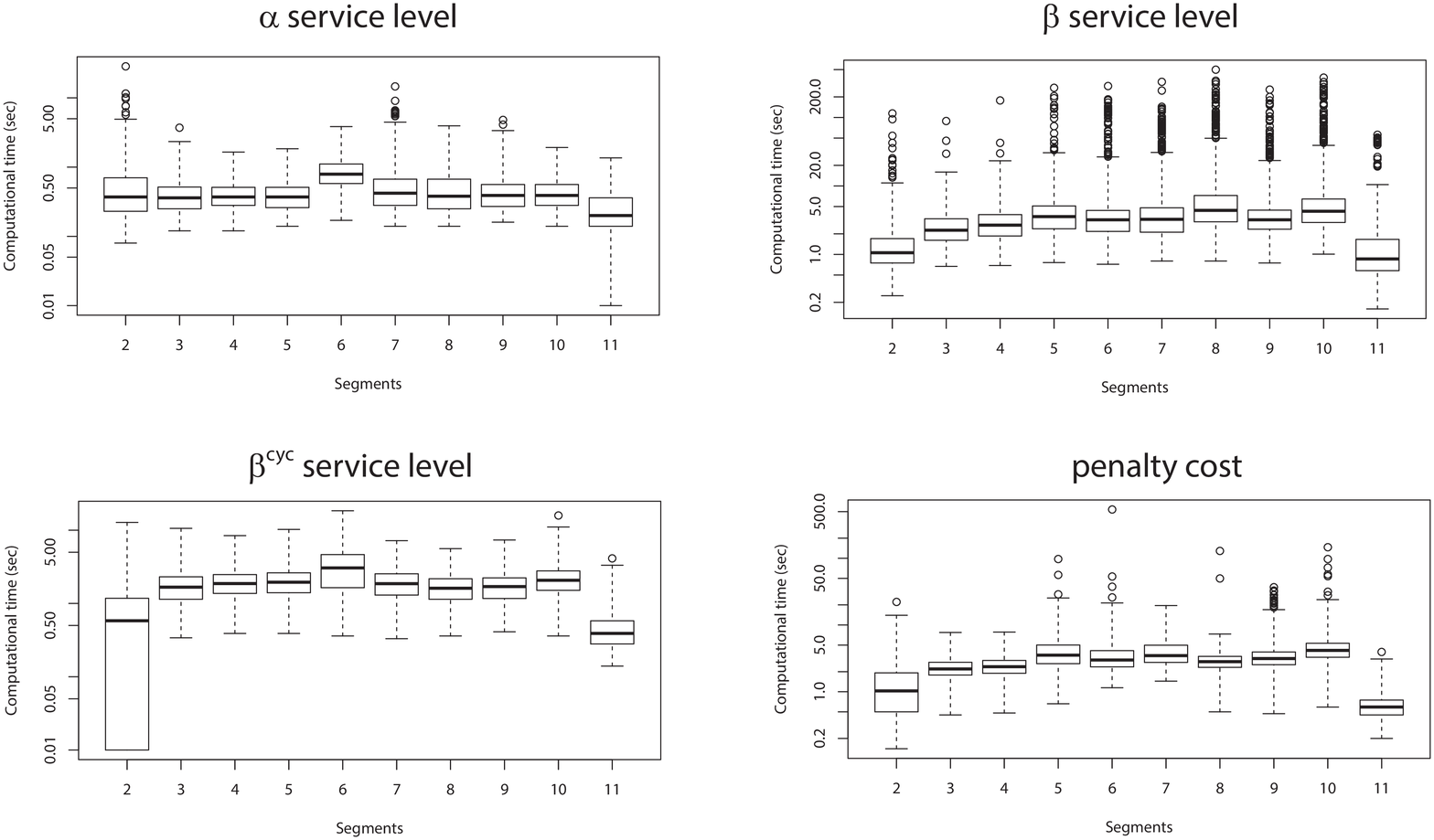}
\caption{Boxplots illustrating the computational time trend for different number of segments used in the piecewise linear approximation. }
\label{fig:cap_6}
\end{figure}

As shown in Fig. \ref{fig:cap_5} also in this case the optimality gap shrinks exponentially fast for all models considered when the number of segments in the linearisation increases. 
Furthermore, as shown in Fig. \ref{fig:cap_6} also in this case all instances could be solved in few seconds.
We do not discuss models that operate under lost sales and generic probability distributions since results obtained were comparable to those already discussed.

\section{Conclusions}\label{sec:conclusions}

We developed MILP formulations for the non-stationary stochastic lot sizing problem. Our formulations exploit a piecewise linearisation of the first order loss function and of its complementary function. We modelled a number of variants of this problem under different service measures: $\alpha$ service level constraints in each period; a penalty cost oriented formulation; a prescribed $\beta^{\text{cyc}}$ service level for each replenishment cycle independently; and a classical $\beta$ service level, as found in the literature. Our models can be easily adapted to operate in a lost sales setting. Our approach has a number of advantages with respect to other existing approaches in the literature. It is versatile, as it enables seamless modelling of several variants of this problem. It is fully linear and, for the special case in which demand in each period is normally distributed, it does not require an offline evaluation of piecewise linearisation coefficients, as these can be derived from a standard table \citep[see e.g.][]{citeulike:12518455}. As shown in our computational experience, viable linearisation parameters for generic demand distributions can be derived by partitioning the support of the demand uniformly or via local search.  Another advantage with respect to other existing approaches is that our models bound from above and below the cost of an optimal plan; by increasing the number of segments in the piecewise linear approximation precision can be improved ad libitum. In our extensive computational experience we demonstrated that optimality gaps shrink exponentially fast in the number of segments used by the piecewise linearisation and that all models developed can be generally solved in a few seconds when up to eleven segments are used in the linearisation.

\subsubsection*{Acknowledgements}
S.A. Tarim and O.A. Kilic were supported by Hacettepe University (HU-BAB) and the ScientiÞc and Technological Research Council of Turkey (TUBITAK) under grant No. 110M500. R. Rossi is supported by the University of Edinburgh CHSS Challenge Investment Fund and by the European Community's Seventh Framework Programme (FP7) under grant no 244994 (project VEG-i-TRADE). 

\bibliographystyle{plainnat}
\bibliography{piecewise}	

\newpage

\section*{Appendix}

In Appendix I, we present a list of all symbols used in the paper; 
in Appendix II, we briefly discuss the two approximate methods adopted in this work for computing good linearisation parameters;
in Appendix III, we provide an overview of the MILP models discussed;
in Appendix IV we discuss special cases that arise when demand is normally distributed.

\subsection*{Appendix I}

\begin{tabular}{p{0.2\columnwidth}p{0.8\columnwidth}}
$x$&a scalar value\\
$\omega$&a random variable\\
$\Omega$&support of  $\omega$\\
$g_{\omega}(x)$&probability density function of $\omega$, where $x\in\Omega$\\
$\tilde{\omega}$&expected value of $\omega$\\
$\sigma_{\omega}$&standard deviation of $\omega$\\
\hline
$Z$&a standard normal random variable\\
$\phi(x)$&standard normal probability density function\\
$\Phi(x)$&standard normal cumulative probability distribution function\\
\hline
$W$&number of regions in a partition of $\Omega$\\
$i$&region index ranging in $1,\ldots,W$\\
$\Omega_i$&a compact region of $\Omega$\\
$p_i$&$\Pr\{\omega\in\Omega_i\}$\\
$\mbox{E}[\omega|\Omega_i]$&conditional expectation of $\omega$ in $\Omega_i$\\
\hline
$\mathcal{L}(x,\omega)$&first order loss function\\
$\widehat{\mathcal{L}}(x,\omega)$&complementary first order loss function\\
$\widehat{\mathcal{L}}_{lb}(x,\omega)$&piecewise linear lower bound of $\widehat{\mathcal{L}}(x,\omega)$ with $W+1$ segments\\
$\widehat{\mathcal{L}}_{ub}(x,\omega)$&piecewise linear upper bound of $\widehat{\mathcal{L}}(x,\omega)$ with $W+1$ segments\\
$e_W$&maximum approximation error for $\widehat{\mathcal{L}}_{lb}(x,\omega)$ and $\widehat{\mathcal{L}}_{ub}(x,\omega)$\\
\hline
$N$&periods in the planning horizon\\
$t$&period index ranging in 1,\ldots,N\\
$d_t$&demand in period $t$\\
$c$&demand coefficient of variation $\tilde{d}_t/\sigma_{d_t}$
\end{tabular}\\
\\
\begin{tabular}{p{0.2\columnwidth}p{0.8\columnwidth}}
$g_t(\cdot)$&probability density function of $d_t$\\
$G_t(\cdot)$&cumulative distribution function of $d_t$\\
$a$&fixed ordering cost, \$$a$ per order\\
$v$&proportional unit cost, \$$v$ per unit ordered\\
$s$&selling price, \$$s$ per item\\
$m$&margin, \$$s-v$ per item\\
$h$&holding cost, \$$h$ per unit of inventory carried to the next period\\
$\mbox{E}[TC]$&expected total cost\\
$I_0$&initial inventory, a scalar\\
$I_t$&inventory level at the end of period $t$\\
$Q_t$&order quantity issued (and received) at the beginning of period $t$\\
$\delta_t$&binary variable set to one if and only if $Q_t>0$.\\
$\alpha$&non stockout probability\\
$\beta$&prescribed fill rate\\
$\beta^{\text{cycle}}$&prescribed cycle fill rate\\
$b$&penalty cost \$$b$ per unit short per unit time\\
$M$&a very large number\\
$P_{jt}$&a binary variable which is set to one if and only if the most recent inventory review before period $t$
was carried out in period $j$\\
$S_t$&order-up-to-level in period $t$\\
$d_{jt}$&random variable representing the convolution $d_j+\ldots+d_t$\\
$G^{-1}_{d_{jt}}(\alpha)$&$\alpha$-quantile of the inverse cumulative distribution function of $d_{jt}$\\
$e_W^{jt}$&maximum approximation error associated with the piecewise linearisation of $\widehat{\mathcal{L}}(x,d_{jt})$\\
$\tilde{I}_t^{lb}$&a lower bound to the true value of $\mbox{E}[\max(I_t,0)]$\\
$\tilde{I}_t^{ub}$&an upper bound to the true value of $\mbox{E}[\max(I_t,0)]$\\
$\tilde{B}_t^{lb}$&a lower bound to the true value of $\mbox{E}[-\min(I_t,0)]$\\
$\tilde{B}_t^{ub}$&an upper bound to the true value of $\mbox{E}[-\min(I_t,0)]$\\
$\tilde{C}_t^{lb}$& a lower bound to the expected total backorders within the replenishment cycle that ends at period t, if any exists. $\tilde{C}_t^{lb}$ is equal to $\tilde{B}_t^{lb}$ if $t$ is the last period of a replenishment cycle\\
$\tilde{C}_t^{ub}$& an upper bound to the expected total backorders within the replenishment cycle that ends at period t, if any exists. $\tilde{C}_t^{ub}$ is equal to $\tilde{B}_t^{ub}$ if $t$ is the last period of a replenishment cycle
\end{tabular}

\subsection*{Appendix II}

We briefly discuss implementation details of the two approximate methods employed in this paper to determine linearisation parameters of the complementary first order loss function for a set of random variables.

The first strategy is quite simple. Consider a set of generic random variables $\omega_1,\ldots,\omega_N$ with complementary first order loss functions $\widehat{\mathcal{L}}(x,\omega_1),\ldots, \widehat{\mathcal{L}}(x,\omega_N)$. 
From \eqref{eq:conditional_expectation} it is clear that, once all $p_i$ are fixed, all $\mbox{E}[\omega|\Omega_k]$ are uniquely determined for each loss function $\widehat{\mathcal{L}}(x,\omega_n)$; this computation can be carried out numerically by using off-the-shelf packages such as Mathematica\footnote{\url{http://www.wolfram.com/mathematica/}} or libraries such as SSJ.\footnote{\url{http://www.iro.umontreal.ca/~simardr/ssj/}} If our aim is to derive a piecewise linearisation with $W+1$ segments, then a possible strategy is to partition the support of each random variable $\omega_n$ into $W$ disjoint compact subregions $\Omega_1,\ldots,\Omega_W$, such that $\Pr\{\omega\in \Omega_i\}=1/W$. Once more, these regions can be determined numerically with one of the aforementioned off-the-shelf packages. By recalling that the complementary first order loss function is convex, it follows that the maximum approximation error will be attained at one of the breakpoints of its piecewise linear approximation. Since there are $W$ breakpoints the maximum approximation error associated with the above partitioning can be easily determined for each $\widehat{\mathcal{L}}(x,\omega_1),\ldots, \widehat{\mathcal{L}}(x,\omega_N)$.

The second strategy is still simple, but slightly more refined than the first one. A pseudocode for this strategy is presented in Fig. \ref{fig:srs_coord_desc}. Instead of simply partitioning the support of random variables $\omega_1,\ldots,\omega_N$ uniformly into compact subregions associated with the same probability mass, we now first generate random partitions via simple random sampling and then we try to improve the best of these partitioning via a coordinate descent approach. More specifically, by observing that a partitioning of the support of $\omega_1,\ldots,\omega_N$ is uniquely determined by a tuple of $W$ probabilities $\langle p_1,\ldots,p_W\rangle$ such that $\sum_{i=1}^Wp_i=1$, we can generate $\mathcal{S}$ such tuples (line \ref{alg:srs}) and then pick tuple $\langle p^*_1,\ldots,p^*_W\rangle$ associated with the minimum maximum approximation error over the set of complementary first order loss functions under scrutiny (line \ref{alg:best_srs}). Finally, we try to reduce this error by performing local moves involving each $p^*_i$ separately, for $i=1,\ldots,W-1$ until no further improvement is possible (line \ref{alg:best_coord_desc}). Note that if we decrement/increment $p^*_i$, since $\sum_{i=1}^Wp_i=1$, then one of the other probabilities must be modified accordingly, this explains why the loops ends at $W-1$; in other words, we have $W-1$ degrees of freedom in the search.
\begin{figure}
\begin{algorithm}[H]
\KwData{
random variables $\omega_1,\ldots,\omega_N$; size $\mathcal{S}$ of the population for simple random sampling; step size $\epsilon$ in the coordinate descent.
}
\KwResult{a tuple $\langle p^*_1,\ldots,p^*_W\rangle$}
generate $\mathcal{S}$ tuples $\langle p^k_1,\ldots,p^k_W\rangle$\;
\For{$k=1$ \KwTo $\mathcal{S}$}{
\lnl{alg:srs} determine the minimum maximum approximation error over $\widehat{\mathcal{L}}(x,\omega_1),\ldots, \widehat{\mathcal{L}}(x,\omega_N)$ under partitioning $\langle p^k_1,\ldots,p^k_W\rangle$\;
\If{partitioning $\langle p^k_1,\ldots,p^k_W\rangle$ ensures the best error so far}{
\lnl{alg:best_srs} record $\langle p^k_1,\ldots,p^k_W\rangle$ as $\langle p^*_1,\ldots,p^*_W\rangle$\;
}
}
\lnl{alg:best_coord_desc} \Repeat{the minimum maximum approximation error cannot be improved}{
\For{$i=1$ \KwTo $W-1$}{
	determine which partitioning between $\langle p^*_1,\ldots,p^*_i+\epsilon,\ldots,p^*_W-\epsilon\rangle$ and $\langle p^*_1,\ldots,p^*_i-\epsilon,\ldots,p^*_W+\epsilon\rangle$ leads to the best improvement for the minimum maximum approximation error over $\widehat{\mathcal{L}}(x,\omega_1),\ldots, \widehat{\mathcal{L}}(x,\omega_N)$\;
	update $\langle p^*_1,\ldots,p^*_i,\ldots,p^*_W\rangle$ to reflect the best current partitioning\;
}
}
\caption{Coordinate descent to determine a good partitioning of the random variable supports}
\end{algorithm}
\label{fig:srs_coord_desc}
\end{figure}

\begin{landscape}
\subsection*{Appendix III}
An overview of the MILP models discussed; numbers refer to the respective equations in the text.\\
\\
\resizebox{1\columnwidth}{!}{
\begin{tabular}{p{0.05\columnwidth}|p{0.15\columnwidth}|p{0.1\columnwidth}p{0.3\columnwidth}|p{0.1\columnwidth}p{0.3\columnwidth}}
\multicolumn{2}{c}{}				&\multicolumn{2}{c}{Lower bound}	&\multicolumn{2}{c}{Upper bound}	\\			
\multicolumn{2}{c}{}				&Objective		&Subject to	&Objective		&Subject to	\\
\cline{2-6}		
\multirow{8}{*}{\rotatebox{90}{\mbox{Backorders}}}&						
$\alpha$ service level			&\eqref{eq:obj_lb}	&\eqref{eq:inventory_conservation_milp} \eqref{eq:reordering_milp} \eqref{eq:service_1} \eqref{eq:service_2} \eqref{eq:service_3} \eqref{eq:binary_p_jt} \eqref{eq:binary_delta_t} \eqref{eq:piecewise_compl_loss_lb}	
							&\eqref{eq:obj_ub}	&\eqref{eq:inventory_conservation_milp} \eqref{eq:reordering_milp} \eqref{eq:service_1} \eqref{eq:service_2} \eqref{eq:service_3} \eqref{eq:binary_p_jt} \eqref{eq:binary_delta_t} \eqref{eq:piecewise_compl_loss_ub}\\
\cline{2-6}							
&penalty cost					&\eqref{eq:obj_pen_lb}	&\eqref{eq:inventory_conservation_milp} \eqref{eq:reordering_milp} \eqref{eq:service_2} \eqref{eq:service_3} \eqref{eq:binary_p_jt} \eqref{eq:binary_delta_t} \eqref{eq:piecewise_compl_loss_lb} \eqref{eq:piecewise_loss_lb}  
							&\eqref{eq:obj_pen_ub}	&\eqref{eq:inventory_conservation_milp} \eqref{eq:reordering_milp} \eqref{eq:service_2} \eqref{eq:service_3} \eqref{eq:binary_p_jt} \eqref{eq:binary_delta_t} \eqref{eq:piecewise_compl_loss_ub} \eqref{eq:piecewise_loss_ub} \\
\cline{2-6}								
&$\beta^{\text{cyc}}$ service level	&\eqref{eq:obj_lb}	&\eqref{eq:inventory_conservation_milp} \eqref{eq:reordering_milp} \eqref{eq:service_2} \eqref{eq:service_3} \eqref{eq:binary_p_jt} \eqref{eq:binary_delta_t} \eqref{eq:piecewise_compl_loss_lb} \eqref{eq:piecewise_loss_lb} \eqref{eq:fill_rate_lb}
							&\eqref{eq:obj_ub}	&\eqref{eq:inventory_conservation_milp} \eqref{eq:reordering_milp} \eqref{eq:service_2} \eqref{eq:service_3} \eqref{eq:binary_p_jt} \eqref{eq:binary_delta_t} \eqref{eq:piecewise_compl_loss_ub} \eqref{eq:piecewise_loss_ub} \eqref{eq:fill_rate_ub}\\
\cline{2-6}								
&$\beta$ service level				&\eqref{eq:obj_lb}	&\eqref{eq:inventory_conservation_milp} \eqref{eq:reordering_milp} \eqref{eq:service_2} \eqref{eq:service_3} \eqref{eq:binary_p_jt} \eqref{eq:binary_delta_t} \eqref{eq:piecewise_compl_loss_lb} \eqref{eq:piecewise_loss_lb} \eqref{eq:end_of_cycle_backorders_lb} \eqref{eq:classical_fill_rate_lb} 
							&\eqref{eq:obj_ub}	&\eqref{eq:inventory_conservation_milp} \eqref{eq:reordering_milp} \eqref{eq:service_2} \eqref{eq:service_3} \eqref{eq:binary_p_jt} \eqref{eq:binary_delta_t} \eqref{eq:piecewise_compl_loss_ub} \eqref{eq:piecewise_loss_ub} \eqref{eq:end_of_cycle_backorders_ub} \eqref{eq:classical_fill_rate_ub}\\
							\hline\hline
\multirow{8}{*}{\rotatebox{90}{\mbox{Lost sales}}}&							
$\alpha$ service level			&\eqref{eq_obj_ls_lb}	&\eqref{eq:piecewise_compl_loss_lb} \eqref{eq:inventory_conservation_milp_ls} \eqref{eq:reordering_milp_ls_1a} \eqref{eq:reordering_milp_ls_1b} \eqref{eq:reordering_milp_ls_3} \eqref{eq:reordering_milp_ls_2} \eqref{eq:service_1_ls} \eqref{eq:service_2_ls} \eqref{eq:service_3_ls} \eqref{eq:binary_p_jt_ls} \eqref{eq:binary_delta_t_ls}
							&					&replace $\tilde{I}^{ub}_t$, \eqref{eq:piecewise_compl_loss_lb} with $\tilde{I}^{lb}_t$, \eqref{eq:piecewise_compl_loss_ub}\\
\cline{2-6}								
&penalty cost					&\eqref{eq_obj_ls_pen_lb}&\eqref{eq:piecewise_compl_loss_lb} \eqref{eq:piecewise_loss_lb}  \eqref{eq:inventory_conservation_milp_ls} \eqref{eq:reordering_milp_ls_1a} \eqref{eq:reordering_milp_ls_1b} \eqref{eq:reordering_milp_ls_3} \eqref{eq:reordering_milp_ls_2} \eqref{eq:service_2_ls} \eqref{eq:service_3_ls} \eqref{eq:binary_p_jt_ls} \eqref{eq:binary_delta_t_ls}
							&					&replace $\tilde{I}^{ub}_t$, \eqref{eq:piecewise_compl_loss_lb}, \eqref{eq:piecewise_loss_lb}, with $\tilde{I}^{lb}_t$, \eqref{eq:piecewise_compl_loss_ub}, \eqref{eq:piecewise_loss_ub}\\
\cline{2-6}	
&$\beta^{\text{cyc}}$ service level	&\eqref{eq_obj_ls_lb}	&\eqref{eq:piecewise_compl_loss_lb} \eqref{eq:piecewise_loss_lb} \eqref{eq:fill_rate_lb} \eqref{eq:inventory_conservation_milp_ls} \eqref{eq:reordering_milp_ls_1a} \eqref{eq:reordering_milp_ls_1b} \eqref{eq:reordering_milp_ls_3} \eqref{eq:reordering_milp_ls_2} \eqref{eq:service_2_ls} \eqref{eq:service_3_ls} \eqref{eq:binary_p_jt_ls} \eqref{eq:binary_delta_t_ls}
							&					&replace $\tilde{I}^{ub}_t$, \eqref{eq:piecewise_compl_loss_lb}, \eqref{eq:piecewise_loss_lb}, \eqref{eq:fill_rate_lb} with $\tilde{I}^{lb}_t$, \eqref{eq:piecewise_compl_loss_ub}, \eqref{eq:piecewise_loss_ub}, \eqref{eq:fill_rate_ub}\\			
\cline{2-6}												
&$\beta$ service level				&\eqref{eq_obj_ls_lb}	&\eqref{eq:piecewise_compl_loss_lb} \eqref{eq:piecewise_loss_lb} \eqref{eq:end_of_cycle_backorders_lb} \eqref{eq:classical_fill_rate_lb} \eqref{eq:inventory_conservation_milp_ls} \eqref{eq:reordering_milp_ls_1a} \eqref{eq:reordering_milp_ls_1b} \eqref{eq:reordering_milp_ls_3} \eqref{eq:reordering_milp_ls_2} \eqref{eq:service_2_ls} \eqref{eq:service_3_ls} \eqref{eq:binary_p_jt_ls} \eqref{eq:binary_delta_t_ls}
							&					&replace $\tilde{I}^{ub}_t$, \eqref{eq:piecewise_compl_loss_lb}, \eqref{eq:piecewise_loss_lb}, \eqref{eq:end_of_cycle_backorders_lb}, \eqref{eq:classical_fill_rate_lb} with $\tilde{I}^{lb}_t$, \eqref{eq:piecewise_compl_loss_ub}, \eqref{eq:piecewise_loss_ub}, \eqref{eq:end_of_cycle_backorders_ub}, \eqref{eq:classical_fill_rate_ub} 					
\end{tabular}}
\end{landscape}

\subsection*{Appendix IV}

In this appendix we discuss model variants for the case in which demand in each period is normally distributed. The assumption that demand is normally distributed plays a prominent role in inventory theory \cite[see e.g.][]{citeulike:10186133} and most of existing works on stochastic lot sizing focus on this distribution. An important property of the first order normal loss function is that it can be derived from its standard counterpart. Let $\omega$ be a normally distributed random variable with mean $\mu$ and standard deviation $\sigma$. Let $\phi(x)$ be the standard normal probability density function and $\Phi(x)$ the respective cumulative distribution function.
\begin{lem}\label{lem:folf_std_norm}
The complementary first order loss function of $\omega$ can be expressed in terms of the standard normal cumulative distribution function as
\begin{equation}\label{eq:compl_fol_norm}
\widehat{\mathcal{L}}(x,\omega)=\sigma\int_{-\infty}^{\frac{x-\mu}{\sigma}} \Phi(t)\,dt=\sigma\widehat{\mathcal{L}}\left(\frac{x-\mu}{\sigma},Z\right),
\end{equation}
where $Z$ is a standard normal random variable.
\end{lem}
\citet{citeulike:12518455} discussed how to obtain an optimal partitioning of the support under a framework that minimises the maximum approximation error. The same work reports standard linearisation parameters for the case in which $\omega$ is a standard normal random variable.

\subsubsection*{$\alpha$ service level constraints}

If demand in each period is normally distributed, we can exploit Lemma \ref{lem:folf_std_norm} to reduce the number of linearisation parameters that are needed in the model. Consider a partition of the support $\Omega$ of a standard normal random variable $Z$ into $W$ adjacent regions $\Omega_i$. Recall that $p_i=\Pr\{Z\in\Omega_i\}$, by exploiting the Jensen's piecewise linear lower bound, we introduce the following constraints in the model
\[\tilde{I}^{lb}_t\geq \tilde{I}_t \sum_{k=1}^i p_k - \sum_{j=1}^t \left( \sum_{k=1}^i p_k\mbox{E}[Z|\Omega_i]\right) P_{jt}\sigma_{d_{jt}}~~~t=1,\ldots,N;~i=1,\ldots,W \]
where $\sigma_{d_{jt}}$ denotes the standard deviation of $d_j+\ldots+d_t$ and $\tilde{I}^{lb}_t\geq0$. This expression follows immediately from Lemma \ref{lem:piecewise_linear_lb} and Lemma \ref{lem:folf_std_norm}: consider a replenishment cycle covering periods $j,\ldots,t$ and associated order-up-to-level $S$. We aim to enforce $\tilde{I}^{lb}_t\geq \sigma \widehat{\mathcal{L}}^i_{lb}\left((S-\tilde{d}_{jt})/\sigma_{d_{jt}},Z\right)$ for all $i=1,\ldots,W$. Observe that $S-\tilde{d}_{jt}=\tilde{I}_t$, the above expression follows immediately. We then derive $\tilde{I}^{ub}_t$  from $\tilde{I}_t$
\[\tilde{I}^{ub}_t\geq \tilde{I}_t \sum_{k=1}^i p_k + \sum_{j=1}^t \left(e_W- \sum_{k=1}^i p_k\mbox{E}[Z|\Omega_i]\right) P_{jt}\sigma_{d_{jt}}~~~
\begin{array}{l}
t=1,\ldots,N,\\
i=1,\ldots,W; 
\end{array}\]
where  $\tilde{I}^{ub}_t\geq \sum_{j=1}^t e_W P_{jt}\sigma_{d_{jt}}$ for $t=1,\ldots,N$; and $e_W$ denotes the maximum approximation error associated with a partition comprising $W$ regions; linearisation parameters that minimise the maximum approximation error $e_W$ for a given number $W+1$ of segments, where $W=1,\ldots,10$, can be found in \cite{citeulike:12518455}.

A final remark that is worth making is the fact that it is possible to replace $e_W$ with $e_W/2$ in the above equations to obtain an approximation of the first order loss function that minimises the maximum absolute error. However, this approximation does not allow one to establish if the cost produced by the model is an upper or a lower bound for the true cost of an optimal plan. 

\subsubsection*{Penalty cost scheme}

We discuss how to handle the case in which demand is normally distributed by exploiting Lemma \ref{thm:relationship_fol_com_1}. 
We obtain $\tilde{B}^{lb}_t$ and $\tilde{B}^{ub}_t$  from $\tilde{I}_t$ by exploiting the connection between the Jensen's piecewise linear lower bound and the piecewise linear upper bound to the first order loss function (Lemma \ref{thm:relationship_fol_com_1}).

\begin{equation}\label{eq:piecewise_loss_lb_norm}
\tilde{B}^{lb}_t\geq - \tilde{I}_t + \tilde{I}_t \sum_{k=1}^i p_k - \sum_{j=1}^t \left( \sum_{k=1}^i p_k\mbox{E}[Z|\Omega_i]\right) P_{jt}\sigma_{d_{jt}}~~~
\begin{array}{l}
t=1,\ldots,N,\\
i=1,\ldots,W; 
\end{array}
\end{equation}
where $\tilde{B}^{ub}_t\geq  - \tilde{I}_t$ and 
\begin{equation}\label{eq:piecewise_loss_ub_norm}
\tilde{B}^{ub}_t\geq -\tilde{I}_t + \tilde{I}_t \sum_{k=1}^i p_k + \sum_{j=1}^t \left(e_W- \sum_{k=1}^i p_k\mbox{E}[Z|\Omega_i]\right) P_{jt}\sigma_{d_{jt}}~~~
\begin{array}{l}
t=1,\ldots,N,\\
i=1,\ldots,W; 
\end{array}
\end{equation}
where  $\tilde{B}^{ub}_t\geq  - \tilde{I}_t + \sum_{j=1}^t e_W P_{jt}\sigma_{d_{jt}}$.  

\end{document}